\documentclass[12pt]{amsart}

\usepackage{enumerate}
\usepackage{amsmath,amsthm,verbatim,amssymb,amsfonts,amscd,graphicx}
\usepackage{graphics}
\usepackage{hyperref}
\usepackage{multicol}
\usepackage{mathrsfs}
\usepackage[left=2cm,right=2cm,top=2.5cm,bottom=2.5cm]{geometry}
\usepackage[latin2]{inputenc}

\usepackage{tikz}
\usepackage{pgfplots}
\pgfplotsset{compat = newest}

\theoremstyle{plain}
\newtheorem{theorem}{Theorem}[section]
\newtheorem{corollary}[theorem]{Corollary}
\newtheorem{lemma}[theorem]{Lemma}
\newtheorem{proposition}[theorem]{Proposition}

\newtheorem*{algorithm*}{Algorithm}
\newtheorem{question}[theorem]{Question}

\theoremstyle{definition}

\newtheorem{s}[theorem]{}
\newtheorem{remark}[theorem]{Remark}

\newtheorem{example}[theorem]{Example}
\newtheorem*{example*}{\it Example}

\theoremstyle{remark}
\newtheorem*{claim*}{\it Claim}

\newtheorem*{case*}{\it Case}
\newtheorem*{note*}{\it Note}

\newcommand{\NN}{\mathbb{N}}

\begin{document}
\title{Computing epsilon multiplicities in graded algebras}
\author{Suprajo Das}
\address{Indian Institute of Technology Bombay, Mumbai, 400076, India}
\email{dassuprajo@gmail.com}
\author{Saipriya Dubey}
\address{Indian Institute of Technology Dharwad, Karnataka, 580011, India}
\email{saipriya721@gmail.com}
\author{Sudeshna Roy}
\address{Tata Institute of Fundamental Research, Mumbai, 400005, India}
\email{sudeshnaroy.11@gmail.com}
\author{Jugal K. Verma}
\address{Indian Institute of Technology Bombay, Mumbai, 400076, India}
\email{verma.jugal@gmail.com}

\begin{abstract}
This article investigates the computational aspects of the $\varepsilon$-multiplicity. Primarily, we show that the $\varepsilon$-multiplicity of a homogeneous ideal $I$ in a two-dimensional standard graded domain of finite type over an algebraically closed field of arbitrary characteristic, is always a rational number. In this situation, we produce a formula for the $\varepsilon$-multiplicity of $I$ in terms of certain mixed multiplicities associated to $I$. In any dimension, under the assumptions that the saturated Rees algebra of $I$ is finitely generated, we give a different expression of the $\varepsilon$-multiplicity in terms of mixed multiplicities by using the Veronese degree. This enabled us to make various explicit computations of $\varepsilon$-multiplicities. We further write a Macaulay2 algorithm to compute $\varepsilon$-multiplicity (under the Noetherian hypotheses) even when the base ring is not necessarily standard graded.
\end{abstract}

\maketitle

\section{Introduction}
The primary aim of this article is to explicitly compute the $\varepsilon$-multiplicities of certain homogeneous ideals $I$ in a finitely generated graded algebra $R$ over a field. This is based on B. Ulrich's suggestion, who introduced this invariant together with J. Validashti. Using standard geometric arguments, we give a complete description when $R$ is a two-dimensional standard graded domain over a field, by means of some well-known invariants that can be computed using Macaulay2. We still have very little understanding beyond dimension two. Nevertheless, we produce a Macaulay2 algorithm to compute the $\varepsilon$-multiplicity under the assumption that the saturated Rees algebra of $I$ is finitely generated. A few examples have been given utilizing our results.

Let $(R,\mathfrak{m})$ be a Noetherian local ring of dimension $d$. Let $\lambda_R(N)$ denote the length of a finite length $R$-module $N$. The Hilbert-Samuel multiplicity of an $\mathfrak{m}$-primary ideal $I$ in $R$ is defined as $$e(I) = \lim_{n\to\infty}\dfrac{\lambda_R(R/I^n)}{n^d/d!} = \lim_{n\to\infty}\dfrac{\lambda_R(I^n/I^{n+1})}{n^{d-1}/(d-1)!}.$$ Classically, the Hilbert-Samuel multiplicity is used to define the intersection numbers for varieties. Other significant applications include Rees' criterion for integral dependence. Being an indispensable technique for the study of local rings, intersection theory and singularity theory, it became necessary to relax the restriction on $I$ to be $\mathfrak{m}$-primary. This led to many new numerical invariants over the last few decades. However, each has its own advantages and disadvantages.

The idea of $\varepsilon$-multiplicity originates in the works of S. Kleiman, B. Ulrich and J. Validashti \cite{KUV}. This invariant has close relations with the volume of divisors and has found applications in the equisingularity theory. Currently, using the theory of integral closure of modules, researchers are trying to obtain algebraic conditions for equisingularity.

Given an ideal $I$ in a $d$-dimensional Noetherian local ring $(R,\mathfrak{m})$, its $\varepsilon$-multiplicity is defined as $$\varepsilon(I) := \limsup_{n \to \infty}\dfrac{\lambda_R\left(H^0_{\mathfrak{m}}\left(R/I^{n}\right)\right)}{n^d/d!} = \limsup_{n \to \infty}\dfrac{\lambda_R\left(\widetilde{I^n}/I^n\right)}{n^d/d!},$$ where $\widetilde{I^n} := I^n\colon_R \mathfrak{m}^{\infty}$ denotes the saturation of the ideal $I^n$, see \cite{BJ} or \cite{BJ2}. S. D. Cutkosky \cite{DC6} showed that `$\limsup$' in the definition can be replaced by `$\lim$' if the ring $R$ is analytically unramified. Note that $\varepsilon(I)$ is bounded above by an integer $j(I)$, the so-called $j$-multiplicity of $I$, defined as $$j(I) = \lim_{n \to \infty}\dfrac{\lambda_R\left(H^0_{\mathfrak{m}}\left(I^n/I^{n+1}\right)\right)}{n^{d-1}/(d-1)!}.$$ The $j$-multiplicity is relatively easier to compute since $\oplus_{n\geq 0}H^0_{\mathfrak{m}}(I^{n}/I^{n+1})t^n$ is a finitely generated module over the Rees algebra $\mathcal{R}(I):=\oplus_{n \geq 0} I^nt^n$ of $I$ and enjoys many properties similar to those of the usual multiplicity, see \cite{NU10}. However, the $\varepsilon$-multiplicity is rather mysterious. In fact, S. D. Cutkosky et al. \cite{DC3} have constructed ideals in a four-dimensional regular local ring for which the $\varepsilon$-multiplicity is an irrational number. This shows that the $\varepsilon$-function need not even have polynomial growth eventually and cannot arise from a finitely generated graded module.

Under the assumptions that the saturated Rees algebra $\mathcal{S}(I):=\oplus_{n\geq 0}\widetilde{I^n}t^n$ of $I$ is Noetherian, J. Herzog, T. J. Puthenpurakal and J. K. Verma \cite{J} showed that the $\varepsilon$-function of $I$, i.e., $$\varepsilon_I(n) := \lambda_R\left(H^0_{\mathfrak{m}}\left(R/I^n\right)\right) = \lambda_R\left(\widetilde{I^n}/I^n\right)$$ is eventually a quasi-polynomial with constant leading coefficient. In this case, the $\varepsilon$-multiplicity is a rational number but it need not always be an integer. If $I$ is a monomial ideal in a polynomial ring over a field then it follows from \cite{HHT} that the algebra $\mathcal{S}(I)$ is Noetherian and consequently its $\varepsilon$-multiplicity must be a rational number. J. Jeffries and J. Monta\~no \cite{JM13} improved this result by describing the $\varepsilon$-multiplicity of a monomial ideal as the normalized volume of a region with rational vertices. They further used this description to compute $\varepsilon$-multiplicities of some monomial ideals. Recently, J. Jeffries, J. Monta\~no and M. Varbaro \cite{JMV15} computed the $\varepsilon$-multiplicity of ideals defining determinantal varieties using the standard monomial theory. However, there is still an evident lack of examples whose $\varepsilon$-multiplicity is known.

In Section \ref{sectmon} of this paper, we extend a result of J. Herzog, T. Hibi and N. V. Trung \cite[Theorem $3.2$]{HHT} to ideals generated by monomials in system of parameters. In Section \ref{epdim1}, we study the $\varepsilon$-function in one-dimensional local rings. In this situation, the $\varepsilon$-multiplicity agrees with the $j$-multiplicity, see Theorem \ref{epsdim1}. For now, we consider the following graded setup: $R=\oplus_{m\geq 0}R_m$ is a standard graded Noetherian algebra over a field $R_0=k$ with unique homogeneous maximal ideal $\mathfrak{m}$ and $I\subseteq R$ is a homogeneous ideal. The $\varepsilon$-multiplicity of $I$ can be analogously defined in this graded situation.

In Section \ref{epdim2}, we explore the behaviour of the $\varepsilon$-function under the additional assumptions that $R$ is a two-dimensional domain with $k=\overline{k}$. Following the ideas in \cite{DC3} and using the standard literature on projective curves, we show that the $\varepsilon$-function is asymptotically given by $$\varepsilon_I(n) = a_0n^2 + a_1(n)n + a_2(n),$$ where $a_0\in \mathbb{Q}_{\geq 0}$ is a constant, $a_1 \colon \mathbb{N} \to \mathbb{Q}$ is a periodic function and $a_2 \colon \mathbb{N} \to \mathbb{Q}$ is a bounded function, see Theorem \ref{rationalep}. Here, the $\varepsilon$-multiplicity is always a rational number. If we further assume that $\mathrm{char}\;k=0$ then the last term $a_2(n)$ is also shown to be periodic, i.e, $\varepsilon_I(n)$ is eventually a quasi-polynomial with constant leading term, but an example due to D. Rees \cite{rees} shows that the algebra $\mathcal{S}(I)$ can be infinitely generated. Also, eventual periodicity of $a_2(n)$ may fail when $\mathrm{char}\;k>0$ and $k$ has positive transcendence degree over its prime subfield. The counterexample to Zariski's Riemann-Roch problem in positive characteristic, due to S. D. Cutkosky and V. Srinivas \cite[Section $6$]{CS}, can be used to construct such an example where $a_2(n)$ is bounded but not eventually periodic.

In higher dimensions, the behaviour of the $\varepsilon$-function can become very complicated, as suggested by Cutkosky's examples \cite[Example $4.4$]{CHT} and \cite[Section $3$]{DC3}. So, in the rest of the article we assume that the algebra $\mathcal{S}(I)$ is finitely generated.

Classically, the mixed multiplicities, denoted $e_i(\mathfrak{m}\vert I)$, of $\mathfrak{m}$ and $I$ are defined as the normalized leading coefficients of the bivariate polynomial associated to the length function $\lambda_R\left(\frac{\mathfrak{m}^mI^n}{\mathfrak{m}^{m+1}I^n}\right)$, see \cite{Tei73}, \cite{KV}. There is another notion of mixed multiplicities in \cite{hoang}, which comes from the natural bigraded structure on the Rees algebra $\mathcal{R}(I)$ of $I$. It is also remarked in \cite{hoang} that these two types of mixed multiplicities can be related if the ideal $I$ is equigenerated. In Section \ref{Mixedmult} of our paper, we obtain a similar relation even when $I$ is generated in different degrees. However, we need to consider the mixed multiplicities of $\mathfrak{m}$ and a truncated ideal $I_{\geq u}$ for some $u$, which is equigenerated. As a consequence, we are able to express the $\varepsilon$-multiplicity of $I$ in terms of mixed multiplicities whenever at least one of the following two conditions hold:
\begin{enumerate}
 \item[$(i)$] $R$ is a two-dimensional graded domain and $k$ is algebraically closed.
 \item[$(ii)$] The algebra $\mathcal{S}(I)$ is Noetherian.
\end{enumerate}
We utilize this description to produce a two-dimensional graded $k$-algebra (Example \ref{exdense}) with $\varepsilon(I)=\frac{p^2}{q}$ for any two previously decided coprime integers $p$ and $q$. In particular, the set of all possible $\varepsilon$-multiplicities form a dense set in $\mathbb{R}_{\geq 0}$. Based on our Macaulay2 computations, we also give a conjectural expression of the $\varepsilon$-multiplicity for certain `fat' points on an elliptic curve. For computational purposes, we make essential use of the package `{\sf MixedMultiplicity}', written by the third and fourth authors along with others \cite{GMRV23}, to compute the relevant mixed multiplicities using Macaulay2.

Nevertheless, the truncation of the ideal increases the number of generators and therefore the ``mixedMultiplicity" command may fail to produce any output. We reserve the last Section \ref{M2_algorithm} to discuss a Macaulay2 algorithm to compute the epsilon function $\varepsilon_I(n)$ and thereby $\varepsilon(I)$. This algorithm works even when the ring $R$ is not necessarily standard graded. However, to use the algorithm one needs to know beforehand a bound on the generating degrees of $\mathcal{S}(I)$. Nevertheless, we could explicitly compute $\varepsilon$-multiplicity for some concrete examples, especially for ideals defining monomial space curves, which are height two quasi-homogeneous prime ideals in a polynomial ring with three variables. One can verify that $\widetilde{I^n}=I^{(n)}$ for every $n \geq 1$ whenever $\mathrm{ht}\; I=d-1$ and $\mathfrak{m}$ is not associated to $I$. Plenty of research has been done surrounding the finite generation of the symbolic Rees algebra $\oplus_{n \geq 0} I^{(n)}t^n$ over the last three decades. In \cite{herzog06}, Herzog found the length of $I^{(n)}/I^n$ in terms of $\mathrm{Ext}$-modules. For a height two prime ideal $I$ in a regular local ring $R$ with specific presentation matrix, he explicitly computed it, which is indeed $\varepsilon_I(n)$, for some values of `$n$' using homological approach. Utilizing our Macaulay2 algorithm, one can extract the values of $\varepsilon_I(n)$ for every $n \geq 0$, when $I$ is associated to a monomial space curve.

\section{Notations and preliminaries}\label{notprel}

In this section, we recall some background material that we need for our study. We also introduce some notations which we use throughout this article.

\s \emph{Epsilon multiplicity}. Suppose that $(R,\mathfrak{m})$ is a Noetherian local ring of dimension $d$ and $I\subseteq R$ is an ideal. The \emph{Rees algebra} of $I$ is the graded ring $\mathcal{R}(I):=R[It]=\oplus_{n\geq 0}I^nt^n$. The \emph{saturated Rees algebra} of $I$ is the graded ring $\mathcal{S}(I):=\oplus_{n\geq 0}\widetilde{I^n}t^n$, where $\widetilde{I^n} := I^n\colon_R \mathfrak{m}^{\infty}$ denotes the \emph{saturation} of the ideal $I^n$. The \emph{$\varepsilon$-multiplicity} of $I$ is defined as $$\varepsilon(I) := \limsup_{n \to \infty}\dfrac{\lambda_R\left(H^0_{\mathfrak{m}}\left(R/I^{n}\right)\right)}{n^d/d!} = \limsup_{n \to \infty}\dfrac{\lambda_R\big(\widetilde{I^n}/I^n\big)}{n^d/d!}.$$ The $\varepsilon$-function of $I$ is given by $\varepsilon_I(n) = \lambda_R\left(H^0_{\mathfrak{m}}\left(R/I^{n}\right)\right)$. The \emph{analytic spread} $\ell_R(I)$ of $I$ is defined to be the Krull dimension of the graded $R/\mathfrak{m}$-algebra $\mathcal{R}(I)/\mathfrak{m}\mathcal{R}(I)$. The following are some well-known facts about $\varepsilon$-multiplicities:
\begin{enumerate}
 \item[$(i)$] $\varepsilon(I)>0$ if and only if $\ell_R(I)=d$, \cite[Theorem $4.7$]{KV}.
 \item[$(ii)$] If $I\subseteq J$ is a reduction then $\varepsilon(I) = \varepsilon(J)$, \cite[Proposition $4.4$]{KV}.
 \item[$(iii)$] If $R$ is analytically unramified then $\varepsilon(I)$ exists as a limit, \cite[Corollary $6.3$]{DC6}.
 \item[$(iv)$] If $\mathcal{S}(I)$ is a Noetherian ring then $\varepsilon(I)$ exists as a limit and is a rational number, \cite[Theorem $3.2$]{J}.
 \item[$(v)$] $\varepsilon(I)$ can be an irrational number, \cite[Theorem $3.2$]{DC3}.
\end{enumerate}

Analogous definitions and results hold in the graded case as well.

\s \emph{Mixed multiplicities}. Suppose that $R=\oplus_{u\geq 0}R_u$ is a standard graded finitely generated algebra over a field $R_0=k$ with homogeneous maximal ideal $\mathfrak{m} = \oplus_{u\geq 1}R_u$. Set $d= \dim \frac{R}{(0\colon_R I^{\infty})}$. Let $I\subseteq R$ be a nonzero homogeneous ideal in $R$.

\begin{theorem}{\cite[Theorem $4.2$]{hoang}}\label{hoang-trung}
Assume that the ideal $I$ is generated in degrees $b_1\leq \cdots \leq b_s$. Then there exist integers $u_0\geq 0$ and $v_0\geq 0$ such that for all $u\geq b_sv + u_0$ and $v\geq v_0$, the bivariate Hilbert function $\dim_k \left(I^v\right)_u$ is equal to a bivariate numerical polynomial $P(u,v)$ of total degree $d-1$. Moreover, if $P(u,v)$ is written in the form $$P(u,v) = \sum_{i=0}^{d-1} \dfrac{e_i\left(R[It]\right)}{i!(d -1-i)!}u^iv^{d-1-i} + \text{\emph{lower degree terms}},$$ then the coefficients $e_i(R[It])$ are integers for all $i=0,\ldots,d-1$ with $e_{d-1}(R[It]) = e(R)$.
\end{theorem}

N. D. Hoang and N. V. Trung \cite{hoang} call the numbers $e_i(R[It])$ the \emph{mixed multiplicities} of $R[It]$.

We now recollect the definition of classical mixed multiplicities of $\mathfrak{m}$ and $I$. P. Bhattacharya \cite{B57} showed that for all $u\gg 0$ and $v\gg 0$, the bivariate Hilbert function $\dim_k \left(\mathfrak{m}^{u}I^v/\mathfrak{m}^{u+1}I^v\right)$ agrees with a bivariate numerical polynomial $Q(u,v)$ of total degree $d-1$. Moreover, we can write
\begin{equation}\label{KV1}
Q(u,v) = \sum_{i=0}^{d-1} \dfrac{e_i\left(\mathfrak{m} \vert I\right)}{i!(d -1-i)!}u^{d-1-i}v^i + \text{lower degree terms},
\end{equation}
where the coefficients $e_i\left(\mathfrak{m} \vert I\right)$ are integers. D. Katz and J. K. Verma \cite{katz-verma} call the numbers $e_i\left(\mathfrak{m} \vert I\right)$ the \emph{mixed multiplicities} of $\mathfrak{m}$ and $I$.

\s {\bf Notations.}

We now set the notations that we will use throughout this article.

\vspace{0.2cm}
$\boldsymbol{e(N)}$ \hspace{1cm} Hilbert-Samuel multiplicity of the module $N$

$\boldsymbol{\varepsilon(I)}$ \hspace{1.25cm} $\varepsilon$-multiplicity of $I$

$\boldsymbol{\varepsilon_I(n)}$ \hspace{1cm} $\lambda_R\left(H^0_{\mathfrak{m}}(R/I^n)\right)$, the $\varepsilon$-function of $I$

$\boldsymbol{e_i\left(\mathcal{R}(I)\right)}$ \hspace{0.3cm} $i^{th}$-mixed multiplicity of $\mathcal{R}(I)$, as defined in \cite{hoang}

$\boldsymbol{e_j\left(\mathfrak{m}\vert I\right)}$ \hspace{0.5cm} $j^{th}$-mixed multiplicity of $\mathfrak{m}$ and $I$,  as defined in \cite{katz-verma}

$\boldsymbol{H^0_{\mathfrak{m}}(N)}$ \hspace{0.65cm} zeroth local cohomology module of $N$ with respect to $\mathfrak{m}$

$\boldsymbol{I}$ \hspace{1.9cm} homogeneous ideal in $R$

$\boldsymbol{\langle I_{\beta} \rangle}$ \hspace{1.5cm} ideal generated by all degree $\beta$ elements in $I$

$\boldsymbol{\widetilde{I^n}}$ \hspace{1.85cm} saturation of $I^n$ with respect to $\mathfrak{m}$, i.e., $\widetilde{I^n}=\bigcup_{i=1}^{\infty}\big(I^n \colon_R \mathfrak{m}^i\big)$

$\boldsymbol{I^{(n)}}$ \hspace{1.6cm} $n^{th}$-symbolic power of $I$, i.e., $I^{(n)} = \bigcap_{P\in\mathrm{Ass}(R/I)}\left(I^nR_P\cap R\right)$

$\boldsymbol{J}$ \hspace{2.05cm} extension of the ideal $I$ in $S$

$\boldsymbol{k}$ \hspace{2.1cm} field

$\boldsymbol{\lambda_R(N)}$ \hspace{1cm} length of a finite length $R$-module $N$

$\boldsymbol{\ell_R(I)}$ \hspace{1.3cm} analytic spread of an ideal $I$ in $R$

$\boldsymbol{\mathfrak{m}}$ \hspace{2cm} unique homogeneous maximal ideal of $R$

$\boldsymbol{N_{\geq b}}$ \hspace{1.6cm} truncated $R$-submodule $\oplus_{n \geq b} N_n$ of $N$

$\boldsymbol{N_{(\mathfrak{p})}}$ \hspace{1.6cm} homogeneous localization of $N$ with respect to a homogeneous prime ideal $\mathfrak{p}$

$\boldsymbol{\mathfrak{n}}$ \hspace{2.2cm} unique homogeneous maximal ideal of $S$

$\boldsymbol{R}$ \hspace{2.1cm} finitely generated graded $k$-algebra

$\boldsymbol{R^{\langle a \rangle}}$ \hspace{1.65cm} $a^{th}$-Veronese subalgebra $\oplus_{n\geq 0} R_{an}$ of $R$

$\boldsymbol{\mathcal{R}(I)}$ \hspace{1.5cm} Rees algebra of $I$, i.e., $\mathcal{R}(I)=\oplus_{n \geq 0}I^nt^n$

$\boldsymbol{\mathcal{R}_s(I)}$ \hspace{1.35cm} symbolic Rees algebra of $I$, i.e., $\mathcal{R}_s(I)=\oplus_{n \geq 0} I^{(n)} t^n$

$\boldsymbol{S}$ \hspace{2.2cm} graded ring $R[U]$, where $U$ is a variable with $\deg U=1$

$\boldsymbol{\mathcal{S}(I)}$ \hspace{1.6cm} saturated Rees algebra of $I$, i.e., $\mathcal{S}(I)=\oplus_{n \geq 0} \widetilde{I^n} t^n$

$\boldsymbol{\mathrm{Sing}(R)}$ \hspace{0.9cm} the singular locus of $R$

$\boldsymbol{V}$ \hspace{2.15cm} projective scheme $\mathrm{Proj}\;R$

\vspace{0.2cm}
We assume that $(R,\mathfrak{m})$ is a local ring only in Sections \ref{sectmon} and \ref{epdim1}. We further assume that $R$ is standard graded in Sections \ref{Mixedmult} and \ref{epdim2}.

\section{Epsilon multiplicity of monomial ideals}\label{sectmon}

A Noetherian local ring $(R,\mathfrak{m})$ is said to be \emph{formally equidimensional} if $\widehat{R}$ is equidimensional.

\begin{proposition}
Suppose that $(R,\mathfrak{m})$ is a formally equidimensional local ring of characteristic $p>0$ with $\dim R = d>0$. Suppose that $x_1,\ldots ,x_d$ is any system of parameters for $R$ and $I$ is an ideal generated by monomials in $x_i$'s. Then the saturated Rees algebra $\mathcal{S}(I)$ is Noetherian. In particular, $\varepsilon(I)$ exists as a limit and is a rational number.
\end{proposition}

\begin{proof}
We may replace $R$ by its completion $\widehat{R}$ since the relevant lengths do not change when we pass to $\widehat{R}$. Also note that $\widehat{R}$ is equidimensional since $R$ is unmixed. Consider the inclusion $$R^{\prime}:=\mathbb{Z}/p\mathbb{Z}[x_1,\ldots ,x_d] \hookrightarrow R.$$ Then $R^{\prime}$ is isomorphic to a polynomial ring in $d$ variables over the field $\mathbb{Z}/p\mathbb{Z}$. Let $\mathfrak{m}^{\prime} = (x_1,\ldots,x_d)$ be the graded maximal ideal of $R^{\prime}$ and let $J\subseteq R^{\prime}$ be the monomial ideal such that $JR=I$. Then
\begin{equation*}
\left(J^n \colon_{R^{\prime}}{\mathfrak{m}^{\prime}}^{\infty}\right)R \subseteq I^n \colon_{R}\mathfrak{m}^{\infty} \subseteq (I^n)^* \colon_R \mathfrak{m}^{\infty} \subseteq \left(\left(J^n \colon_{R^{\prime}}{\mathfrak{m}^{\prime}}^{\infty}\right)R\right)^*\subseteq \overline{\left(J^n \colon_{R^{\prime}}{\mathfrak{m}^{\prime}}^{\infty}\right)R},
\end{equation*}
where the third inclusion follows from ``colon capturing'' \cite[Theorem $2.3$]{AHS}. In other words, there are graded inclusions of graded $R$-algebras $$A:= \bigoplus\limits_{n\geq 0}\left(J^n \colon_{R^{\prime}}{\mathfrak{m}^{\prime}}^{\infty}\right)Rt^n \subseteq B:= \bigoplus\limits_{n\geq 0}\left(I^n \colon_{R}\mathfrak{m}^{\infty}\right)t^n \subseteq C:= \bigoplus\limits_{n\geq 0} \overline{\left(J^n \colon_{R^{\prime}}{\mathfrak{m}^{\prime}}^{\infty}\right)R}t^n.$$ Since $J$ is a monomial ideal in the polynomial ring $R^{\prime}$ so $A$ is a finitely generated graded $R$-algebra by \cite[Theorem $3.2$]{HHT}. Note that $C\subseteq \overline{A}$, where $\overline{A}$ is the integral closure of $A$ inside $R[t]$. As $R$ is excellent so $\overline{A}$ is a finitely generated graded $R$-algebra and also a finite module over $A$. Therefore, by \cite[Proposition $7.8$]{atiyah}, $C$ is a finitely generated graded $R$-algebra and so is $B$. Then the length function $\lambda_R(\widetilde{I^n}/I^n)$ is eventually given by a numerical quasi-polynomial in $n$ of degree at most $d$ and having constant leading coefficient, see \cite[Theorem $2.4$]{J}. In particular, $\varepsilon(I)$ exists as a limit and $\varepsilon(I)\in\mathbb{Q}$.
\end{proof}

\section{Epsilon function in dimension one}\label{epdim1}

Let $R$ be a Noetherian ring and $I\subseteq R$ an ideal.
The \emph{$n$-th symbolic power} of $I$ is the ideal $$I^{(n)} = \bigcap_{P\in\mathrm{Ass}(R/I)}\varphi^{-1}\left(I^nR_P\right),$$ where $\varphi_P \colon R \to R_P$ is the localization map.

\begin{lemma}\label{stupid1}
Let $(R,\mathfrak{m})$ be a Noetherian local ring and $I\subseteq R$ an ideal such that $I$ has no embedded primes and $\mathrm{ht}(I)=0$. Then there exists an integer $n_0\geq 0$ such that $I^{(n)} = I^{(n+1)}$ for all $n\geq n_0$.
\end{lemma}
\begin{proof}
Let $P_1,\ldots,P_s$ be the minimal primes of $R$, so each $R_{P_i}$ is an Artinian local ring. Thus there exists an integer $n_0\geq 0$ such that $P_i^nR_{P_i} = (0)$ for all $n\geq n_0$ and $i=1,\ldots,s$. Consider the natural maps $\varphi_i \colon R \to R_{P_i}$. By relabelling the $P_i$'s, we may assume that $P_1,\ldots,P_r$ are the minimal primes of $I$, where $r$ is some integer with $1\leq r\leq s$. Then for all $n\geq n_0$, $$I^{(n)} = \bigcap_{i=1}^r \varphi_i^{-1}\left(I^nR_{P_i}\right) = \bigcap_{i=1}^r\mathrm{ker} \;\varphi_i.$$
\end{proof}

\begin{lemma}\label{stupid2}
Let $(R,\mathfrak{m})$ be a Noetherian local ring and $I\subseteq R$ an ideal. Then for all integers $n\geq 0$, $$\widetilde{I^n} = \widetilde{\left(\widetilde{I}\right)^n}.$$
\end{lemma}
\begin{proof}
Clearly $I \subseteq \widetilde{I}$ which implies $I^n \subseteq \left(\widetilde{I}\right)^n$ and thus $\widetilde{I^n} \subseteq \widetilde{\left(\widetilde{I}\right)^n}$. As $\{\widetilde{I^n}\}_{n\in\mathbb{N}}$ is a filtration, so $\left(\widetilde{I}\right)^n \subseteq \widetilde{I^n}$ and hence $\widetilde{\left(\widetilde{I}\right)^n} \subseteq \widetilde{\left(\widetilde{I^n}\right)} = \widetilde{I^n}$.
\end{proof}

\begin{theorem}\label{epsdim1}
Let $(R,\mathfrak{m})$ be a one-dimensional Noetherian local ring and $I\subseteq R$ an ideal. Then there exist integers $a$ and $b$ such that for all integers $n\gg 0$, $$\lambda_R\left(H^0_{\mathfrak{m}}\left(R/I^n\right)\right) = an + b.$$
\end{theorem}
\begin{proof}
The assertion is clear if $\mathrm{ht}(I)=1$, i.e. $I$ is $\mathfrak{m}$-primary. We may thus assume that $\mathrm{ht}(I)=0$. Let $J=\widetilde{I}$. By construction $\mathrm{ht}(J)=0$ and $\mathrm{Min}(R/J) = \mathrm{Ass}(R/J)$. By Lemma \ref{stupid2} $$J^{(n)} = \widetilde{J^n} = \widetilde{I^n}$$ for all $n\geq 0$. We know from Lemma \ref{stupid1} that there exists an integer $n_0\geq 0$ such that $J^{(n)} = J^{(n+1)}$ for all $n\geq n_0$ and denote this stable ideal by $\mathfrak{J}$. Then for all $n\geq n_0$, we have $$H^0_{\mathfrak{m}}\left(R/I^n\right) = \widetilde{I^n}/I^n = \mathfrak{J}/I^n.$$ The inclusion $I^n/I^{n+1} \subseteq \widetilde{I^n}/I^{n+1} = \mathfrak{J}/I^{n+1}$ implies that $I^n/I^{n+1}$ has finite length for all $n\geq n_0$. Since $\dim R=1$, the length function $\lambda_R(I^n/I^{n+1})$ eventually becomes a constant, which equals the $j$-multiplicity of $I$. By increasing $n_0$ if needed, we have for all $n>n_0$, $$\lambda_R\left(H^0_{\mathfrak{m}}\left(R/I^n\right)\right) = \lambda_R\left(\mathfrak{J}/I^n\right) = \lambda_R\left(\mathfrak{J}/I^{n_0}\right) +\sum_{i=n_0+1}^n \lambda_R(I^i/I^{i+1}) =  a(n-n_0) + \lambda_R\left(\mathfrak{J}/I^{n_0}\right) = an+b.$$
\end{proof}

\begin{example}
Let $k$ be a field and $d\geq 2$ an integer. Consider the face ring $R=k[X_1,\ldots,X_d]/I_{\Delta}$, where $I_{\Delta} = \left(X_iX_j \;\colon\; 1\leq i<j\leq d\right)$ is the face ideal of the simplicial complex $\Delta=\{\emptyset,\{1\},\ldots,\{d\}\}$. For each $i=1,\ldots,d$, consider the prime ideal $P_i = (x_1,\ldots,x_{i-1},x_{i+1},\ldots,x_d)$ and let $\mathfrak{m}=(x_1,\ldots,x_d)$. Then for all integers $n\geq 1$, $$\lambda_R\left(H^0_{\mathfrak{m}}\left(R/P_i^n\right)\right) = (d-1)(n-1).$$
\end{example}

\begin{proof}
By symmetry it suffices to prove the statement for $P_1$. We first claim that $\widetilde{P_1^n} = P_1^{(n)} = P_1$ for all $n\geq 1$. To see this, let $j$ be an integer satisfying $2\leq j\leq d$. We have $x_1 \not\in P_1 = \mathrm{rad}\big(P_1^{(n)}\big)$ but $0=x_1x_j \in P_1^{(n)}$ which means that $x_j \in P_1^{(n)}$ as $P_1^{(n)}$ is $P_1$-primary. Since $P_1 \subseteq P_1^{(n)} \subseteq \mathrm{rad}\big(P_1^{(n)}\big) = P_1$ our claim follows. Now, we observe that
\begin{equation*}
\lambda_R\left(H^0_{\mathfrak{m}}\left(R/P_1^n\right)\right) = \lambda_R(P_1/P_1^n) = \dim_k \dfrac{(X_2,\ldots,X_d)}{(X_2^n,\ldots,X_d^n)+I_{\Delta}} = \sum_{j=1}^{n-1}\dim_k \dfrac{(X_2^j,\ldots,X_d^j)+I_{\Delta}}{(X_2^{j+1},\ldots,X_d^{j+1})+I_{\Delta}}.
\end{equation*}
Moreover, it follows from \cite[Corollary $1.1.4$]{HH} that for each $j=1,\ldots,n-1$,
\begin{multline*}
 \dim_k \frac{(X_2^j,\ldots,X_d^j)+I_{\Delta}}{(X_2^{j+1},\ldots,X_d^{j+1})+I_{\Delta}} = \#\big\{(i_1,\ldots,i_d)\in\mathbb{N}^d \mid X_1^{i_1}\cdots X_d^{i_d} \in (X_2^j,\ldots,X_d^j)+I_{\Delta}\\\text{and}\;\;X_1^{i_1}\cdots X_d^{i_d} \not\in (X_2^{j+1},\ldots,X_d^{j+1})+I_{\Delta} \big\}
\end{multline*}
which equals $d-1$. So, we are done.
\end{proof}

\section{Hilbert polynomials of bigraded Rees algebras and mixed multiplicities}\label{Mixedmult}

\s{\bf Setup}.\label{setup} Suppose that $R=\oplus_{u\geq 0}R_u$ is a standard graded finitely generated algebra over a field $R_0=k$ with homogeneous maximal ideal $\mathfrak{m} = \oplus_{u\geq 1}R_u$. Further assume that $R$ is a domain of Krull dimension $d>0$. Let $I\subseteq R$ be a nonzero homogeneous ideal.

The assumption that $R$ is a domain, is not necessary. Most of the results go through by taking $d= \dim \frac{R}{(0\colon_R I^{\infty})}$.

Since $I$ is homogeneous, the Rees algebra $R[It]$ of $I$ has a natural bigraded structure by setting $$R[It]_{(u,v)} = \left(I^v\right)_ut^v$$ for all $(u,v)\in\mathbb{N}^2$. We assume that $I$ is minimally generated by homogeneous elements $f_1,\ldots,f_s$ of degrees $b_1,\ldots,b_s$ respectively with $b_1\leq \cdots \leq b_s$. Nakayama's lemma shows that the sequence $\{b_1,\ldots,b_s\}$ is independent of the choice of a minimal homogeneous generating set of $I$. Clearly, then the (non-standard) bigraded Rees algebra $R[It]$ is generated in bidegrees $\{(1,0),(b_1, 1), \ldots, (b_s, 1)\}$. However, if $I$ is generated in equal degrees, say $b$, then we can give a standard bigraded structure on $R[It]$ by setting
\begin{equation}\label{standardrees}
 R[It]_{(u,v)}={(I^v)}_{bv+u} \qquad \mbox{for every $(u,v) \in \mathbb{N}^2$}.
\end{equation}

If $c>b_s$ is an integer then it immediately follows from Theorem \ref{hoang-trung} that for all $v\gg 0$,
\begin{equation}\label{HT1}
\dim_k \left(I^v\right)_{cv} = \left[\sum_{i=0}^{d-1} \dfrac{e_i\left(R[It]\right)}{i!(d -1-i)!}c^i\right]v^{d-1} + \text{lower degree terms}.
\end{equation}

We now aim to relate the two mixed multiplicities $e_i(R[It])$ and $e_j(\mathfrak{m}\vert I)$ mentioned in Section \ref{notprel}. Observe that the first set of invariants are associated to a non-standard bigraded algebra $R[It]$, whereas, the second set of invariants are associated to a standard bigraded algebra $\mathrm{gr}_M R[It]$.

\begin{lemma}[With hypotheses as in Setup \ref{setup}]\label{Res1}
Assume that the ideal $I$ is generated in equal degrees $b\geq 1$. Then there exist $u_0\geq 0$ and $v_0\geq 0$ such that for all $u\geq u_0$ and $v\geq v_0$ we have $$\dim_k \left(I^v\right)_{bv+u} = \sum_{i=0}^{d-1} \dfrac{e_i\left(\mathfrak{m} \vert I\right)}{i!(d -1-i)!}u^{d-1-i}v^i + \text{\emph{lower degree terms}},$$
\end{lemma}
\begin{proof}
Let $M = \left(\mathfrak{m},It\right)$ be the bigraded maximal ideal of $R[It]$. For each $v\geq 1$, $$M^v = (\mathfrak{m}, It)^v=\mathfrak{m}^v \oplus \mathfrak{m}^{v-1}It \oplus \cdots \oplus \mathfrak{m} I^{v-1}t^{v-1} \oplus \left(\oplus_{j\geq v}I^jt^j\right).$$ So the associated graded ring $\mathrm{gr}_M R[It]:= \oplus_{v\geq 0} M^v/M^{v+1}$ has a natural standard bigrading with $$\left(\mathrm{gr}_M R[It]\right)_{(u,v)} = \dfrac{\mathfrak{m}^{u}I^v}{\mathfrak{m}^{u+1}I^v}$$ for all $(u,v)\in\mathbb{N}^2$. Since $I$ is generated in a single degree $b$, for every $v \geq 1$ we have $I^v={(I^v)}_{\geq bv}$. As $R$ is standard graded so $\mathfrak{m}^u = R_{\geq u}$ for all $u\geq 1$. Thus for any fixed pair $(u,v)\in\mathbb{N}^2$,
\[\mathfrak{m}^u I^v = R_{\geq u} \cdot {(I^v)}_{\geq bv}={(I^v)}_{\geq bv +u}.\]
Therefore,
\[\dim_k \mathrm{gr}_M(R[It])_{(u,v)}=\dim_k\frac{\mathfrak{m}^u I^v}{\mathfrak{m}^{u+1} I^v} = \dim_k \frac{{(I^v)}_{\geq bv +u}}{{(I^v)}_{\geq bv +u+1}} =\dim_k {(I^v)}_{bv +u}.\] The assertion of the lemma is now immediate from \eqref{KV1}.
\end{proof}

It was remarked in \cite[page $326$]{hoang} that if the ideal $I$ is generated in equal degrees then each $e_{i}(\mathfrak{m}\vert I)$ can be expressed in terms of $e_0(R[It]), \ldots, e_{d-1}(R[It])$. Here, we provide a proof of their remark.

\begin{lemma}[With hypotheses as in Setup \ref{setup}]\label{mixedrel}
Assume that the ideal $I$ is generated in equal degrees $b\geq 1$. Then $$e_{i}(\mathfrak{m}\vert I) = \sum_{j=0}^i \binom{i}{j}b^j e_{d-1-i+j}\left(R[It]\right)$$ for each $i=0,\ldots,d-1$.
\end{lemma}
\begin{proof}
Let $c>0$ be an integer. From the expression in Lemma \ref{Res1} it follows that
\begin{equation}\label{local1}
 \lim_{n\to\infty}\dfrac{\dim_k\left(I^v\right)_{bv+cv}}{n^{d-1}/(d-1)!} = \sum_{i=0}^{d-1} \binom{d-1}{i} e_i\left(\mathfrak{m} \vert I\right) c^{d-1-i}.
\end{equation}
Again using the expression in Theorem \ref{hoang-trung} we see that
\begin{align}\label{local2}
 \lim_{n\to\infty}\dfrac{\dim_k\left(I^v\right)_{(b+c)v}}{n^{d-1}/(d-1)!} &= \sum_{j=0}^{d-1}\binom{d-1}{j}e_j(R[It])(c+b)^j\nonumber\\
 &= \sum_{i=0}^{d-1}\left(\sum_{j=d-1-i}^{d-1}\binom{d-1}{j}\binom{j}{d-1-i}b^{j-(d-i-1)}e_j(R[It])\right)c^{d-i-1}\nonumber\\
&= \sum_{i=0}^{d-1}\left(\sum_{j=d-1-i}^{d-1}\binom{d-1}{i}\binom{i}{d-1-j}b^{j-(d-i-1)}e_j(R[It])\right)c^{d-i-1}\nonumber\\
 &= \sum_{i=0}^{d-1}\left(\sum_{j=d-1-i}^{d-1}\binom{i}{d-1-j}b^{j-(d-i-1)}e_j(R[It])\right)c^{d-i-1}\binom{d-1}{i}\nonumber\\
 &= \sum_{i=0}^{d-1}\left(\sum_{k=0}^{i}\binom{i}{k}b^{k}e_{d-1-i+k}(R[It])\right)c^{d-i-1}\binom{d-1}{i}.
\end{align}
By comparing the expressions in $\eqref{local1}$ and $\eqref{local2}$ we get the required expression for $e_{i}(\mathfrak{m}\vert I)$.
\end{proof}

Define $S=R[U]$, where $U$ is a variable with $\deg U=1$. Let $J=IS$ be the extended ideal in $S$ and $\mathfrak{n}=\mathfrak{m}+(U)$ be the homogeneous maximal ideal of $S$. Given an integer $\beta\geq 0$, we denote by $\langle J_{\beta}\rangle$, the ideal in $S$ generated by homogeneous elements in $J$ of degree $\beta$. Note that if $\beta\geq b_s$, then $\langle J_{\beta}\rangle$ equals the truncated ideal $J_{\geq \beta} = J\cap \mathfrak{n}^{\beta}$.

\begin{proposition}[With hypotheses as in Setup \ref{setup}]\label{mixedrel0}
Let $\beta>0$ be an integer such that $I$ is generated in degrees $\leq \beta$. For any integer $c>\beta$ there exists an integer $v_0\geq 0$ (depending on $c$) such that for all integers $v\geq v_0$, $$\sum_{u=0}^{cv}\dim_k \left(I^v\right)_u = \left[\sum_{i=0}^{d} \dfrac{e_i\left(\mathfrak{n} \vert \langle J_{\beta}\rangle\right)}{i!(d-i)!}(c-\beta)^{d-i}\right]v^d + \text{\emph{lower degree terms}}.$$
\end{proposition}
\begin{proof}
Since $c>\beta$, it follows from an argument used in \cite[Theorem $2.3$, page $1185$]{DC3} that for all integer $v\geq 0$, $$\sum_{u=0}^{cv}\dim_k \left(I^v\right)_u = \dim_k \left(J^v\right)_{cv}.$$ Consider the diagonal subalgebra $S[Jt]_{\Delta}$ of $S[Jt]$ along the diagonal $\Delta = \{(cv,v) \mid v\geq 0\}$, i.e., $$S[Jt]_{\Delta} = \oplus_{v\geq 0} \left(J^v\right)_{cv}.$$ As $c>\beta$, so by \cite[Lemma $2.2$]{hoang}, $S[Jt]_{\Delta}$ is a standard graded finitely generated $k$-algebra of dimension $d$. Thus, the Hilbert function $\dim_k \left(J^v\right)_{cv}$ is given by a numerical polynomial in $v$ of degree $d$ for all $v\gg 0$.
\begin{claim*}
The $k$-vector spaces $J_{\alpha}$ and $\langle J_{\beta}\rangle_{\alpha}$ are equal for all integers $\alpha\geq \beta$.
\end{claim*}
Clearly $\langle J_{\beta}\rangle_{\alpha} \subseteq J_{\alpha}$. Now pick a homogeneous element $f\in J_{\alpha}$ where $\alpha \geq \beta$. Write $f = \sum_{i=1}^s a_if_i$ where $\{f_1,\ldots,f_s\}$ is a minimal homogeneous generating set of $J$ with $\deg f_i = b_i$ and $a_1,\ldots,a_s$ are homogeneous elements in $R$ with $\deg a_i = \alpha - b_i$. Since $R$ is standard graded, each $a_i$ may be further written as $a_i = \sum_{j=1}^{r_i} g_{i,j}h_{i,j}$, where $g_{i,j} \in R_{\alpha-\beta}$ and $h_{i,j} \in R_{\beta-b_i}$. From this expression of $a_i$ we have $$f = \sum_{i=1}^s \sum_{j=1}^{r_i} g_{i,j}h_{i,j}f_i,$$ where $h_{i,j}f_i \in J_{\beta}$. In other words, $f\in \langle J_{\beta}\rangle_{\alpha}$ and the proof of our claim is complete.

\vspace{0.2cm}
As $c\geq \beta$, the above claim implies that $${S[\langle J_{\beta}\rangle t]}_{\Delta} = {S[Jt]}_{\Delta}.$$ Now using Lemma \ref{Res1} we see that for all $v\gg 0$,
\begin{align*}
 \dim_k \left(J^v\right)_{cv} = \dim_k {\left(\langle J_{\beta}\rangle^v\right)}_{cv} &= \dim_k {\left(\langle J_{\beta}\rangle^v\right)}_{\beta v+(c-\beta)v}\\
 &= \left[\sum_{i=0}^{d} \dfrac{e_i\left(\mathfrak{n} \vert \langle J_{\beta}\rangle\right)}{i!(d-i)!}(c-\beta)^{d-i}\right]v^d + \text{lower degree terms}.
\end{align*}
\end{proof}

\begin{corollary}[With hypotheses as in Setup \ref{setup}]\label{superficial0}
Suppose that the ideal $I$ is generated in a single degree $b>0$ and fix an integer $c>b$. Then for all $v\gg 0$, $$\sum_{u=0}^{cv}\dim_k \left(I^v\right)_u = \left[\sum_{i=0}^{d-1} \dfrac{e_i\left(\mathfrak{m} \vert I\right)}{i!(d-i)!}(c-b)^{d-i}\right]v^d + \text{\emph{lower degree terms}}.$$
\end{corollary}
\begin{proof}
Here, we use the identity obtained in Proposition \ref{mixedrel0}. Since the ideal $I$ is generated in equal degrees so $\left<J_b\right> = J$. The analytic spread $\ell_S(J)$ of $J$ equals the analytic spread $\ell_R(I)$ of $I$. In other words, $\ell_S(J)\leq d$ and therefore $e_d\left(\mathfrak{n}\vert J\right) = 0$ by \cite[Theorem $2.7$]{katz-verma}. Further, note that $U \in \mathfrak{n} \backslash \mathfrak{n}^2$ is a nonzerodivisor on $S$ and hence superficial for $\mathfrak{n}, J$. Also, the images of the ideals $J=IS$ and $\mathfrak{n}=\mathfrak{m}+(U)$ in $S/(U) \cong R$ are $I$ and $\mathfrak{m}$ respectively. Now using similar arguments as in \cite[Theorem $17.4.6$]{HS} for non $\mathfrak{n}$-primary ideal $J$, we get
\[e_i(\mathfrak{n}, J)=e_S(\mathfrak{n}^{[d+1-i]}, J^{[i]}) = e_R(\mathfrak{m}^{[d-i]}, I^{[i]})=e_i(\mathfrak{m}, I).\]
for all $i=0,\ldots, d-1$.
\end{proof}

\begin{lemma}[With hypotheses as in Setup \ref{setup}]\label{noetheps}
Suppose that $I$ and $I^{\prime}$ are nonzero homogeneous ideals in $R$ such that $I\subseteq I^{\prime}$ and $\lambda_R(I^{\prime}/I)<\infty$. Let $\beta>0$ be an integer such that both $I$ and $I^{\prime}$ are generated in degrees $\leq \beta$. Let $J=IS$ and $J^{\prime}=I^{\prime}S$ be the extended ideals in $S$. Then $$\lim\limits_{v\to\infty}\dfrac{\lambda_R({I^{\prime}}^v/I^v)}{v^d/d!} = e_d\left(\mathfrak{n} \vert \langle J^{\prime}_{\beta}\rangle\right) - e_d\left(\mathfrak{n} \vert \langle J_{\beta}\rangle\right).$$
\end{lemma}
\begin{proof}
By a result of Amao \cite[Theorem $3.2$]{JO}, the length function $\lambda_R({I^{\prime}}^v/I^v)$ eventually agrees with a numerical polynomial in $v$ of degree at most $d$. Since $\lambda_R(I^{\prime}/I)<\infty$, so there exists an integer $c>\beta$ such that $\left({I^{\prime}}^v\right)_u = \left(I^v\right)_u$ for all $u\geq cv$ and $v\geq 0$. Therefore, $$\lim\limits_{v\to\infty}\dfrac{\lambda_R({I^{\prime}}^v/I^v)}{v^d/d!} = \lim\limits_{v\to\infty}\dfrac{\sum_{u=0}^{cv}\dim_k \left({I^{\prime}}^v\right)_u}{v^d/d!} - \lim\limits_{v\to\infty}\dfrac{\sum_{u=0}^{cv}\dim_k \left(I^v\right)_u}{v^d/d!}.$$ Using Proposition \ref{mixedrel0} we see that
\begin{align}\label{amaorel}
 \lim\limits_{v\to\infty}\dfrac{\lambda_R({I^{\prime}}^v/I^v)}{v^d/d!} &= \sum_{i=0}^{d} \binom{d}{i} e_i\left(\mathfrak{n} \vert \langle J^{\prime}_{\beta}\rangle\right)(c-\beta)^{d-i} - \sum_{i=0}^{d} \binom{d}{i} e_i\left(\mathfrak{n} \vert \langle J_{\beta}\rangle\right)(c-\beta)^{d-i}\nonumber\\
 &= \sum_{i=0}^{d} \binom{d}{i} \left[e_i\left(\mathfrak{n} \vert \langle J^{\prime}_{\beta}\rangle\right) - e_i\left(\mathfrak{n} \vert \langle J_{\beta}\rangle\right)\right] (c-\beta)^{d-i}
\end{align}
and this relation holds true for all integers $c\gg 0$. In other words, the coefficient of $c^i$ must be zero for all $i=1,\ldots,d$, since this expression must be independent of $c$. This gives that $$e_i\left(\mathfrak{n} \vert \langle J^{\prime}_{\beta}\rangle\right) = e_i\left(\mathfrak{n} \vert \langle J_{\beta}\rangle\right)$$ for all $i=0,\ldots,d-1,$ and our assertion is now immediate.
\end{proof}

For any positive integer $v_0$ we denote by $\mathcal{S}(I)^{\langle v_0\rangle} = \oplus_{v\geq 0}\widetilde{I^{vv_0}}$ the $v_0$-th Veronese subalgebra of the saturated Rees algebra $\mathcal{S}(I)$ of $I$.

\begin{theorem}[With hypotheses as in Setup \ref{setup}]\label{epNoeth}
Suppose that there exists an integer $v_0>0$ such that $\mathcal{S}(I)^{\langle v_0\rangle}$ is a standard graded $R$-algebra. Let $\beta_{v_0}>0$ be an integer such that both $I^{v_0}$ and $\widetilde{I^{v_0}}$ are generated in degrees $\leq \beta_{v_0}$. Then $$\varepsilon(I) = \lim_{v\to\infty}\dfrac{\lambda_R\left(H^0_{\mathfrak{m}}\left(R/I^v\right)\right)}{v^d/d!} = \dfrac{1}{v_0^d}\cdot \left[e_d\left(\mathfrak{n} \vert \langle \left(\widetilde{I^{v_0}}S\right)_{\beta_{v_0}}\rangle\right) - e_d\left(\mathfrak{n} \vert \langle \left(I^{v_0}S\right)_{\beta_{v_0}}\rangle\right)\right].$$
\end{theorem}
\begin{proof}
The existence of the desired limit follows from \cite[Theorem $2.4$]{J} since $\left\{\widetilde{I^v}\right\}_{v\in\mathbb{N}}$ is a Noetherian filtration. Thus, $$\varepsilon(I) = \lim_{v\to\infty}\dfrac{\lambda_R\left(H^0_{\mathfrak{m}}\left(R/I^{vv_0}\right)\right)}{(vv_0)^d/d!} = \dfrac{1}{v_0^d}\cdot\lim_{v\to\infty}\dfrac{\lambda_R\left(\widetilde{I^{vv_0}}/I^{vv_0}\right)}{v^d/d!} = \dfrac{1}{v_0^d}\cdot\lim_{v\to\infty}\dfrac{\lambda_R\left(\left(\widetilde{I^{v_0}}\right)^v/\left(I^{v_0}\right)^v\right)}{v^d/d!}.$$ Our claim is now immediate from Lemma \ref{noetheps}.
\end{proof}

\begin{example}\label{Tess0}
Let $k$ be a field and consider the monomial ideal $I=(XY,YZ,ZX)$ in the polynomial ring $R=k[X,Y,Z]$ with $\mathfrak{m}=(X,Y,Z)$. It follows from \cite[Proposition $5.3$ and Example $4.7$]{HHT} that for $v>0$, $$\widetilde{I^{2v}} = \left(\widetilde{I^2}\right)^v = \left(I^2 + (XYZ)\right)^v = \left(XYZ, X^2Y^2, Y^2Z^2, X^2Z^2\right).$$ Then $I^2$ and $\widetilde{I^2}$ are generated in degrees $\{4\}$ and $\{3,4\}$ respectively. From Theorem \ref{epNoeth} we have ($v_0=2$ and $\beta_{v_0}=4$) $$\varepsilon(I) = \dfrac{1}{2^3}\cdot \left[e_3\left(\mathfrak{n} \vert \langle \left(\widetilde{I^{2}}S\right)_{4}\rangle\right) - e_3\left(\mathfrak{n} \vert \langle \left(I^{2}S\right)_{4}\rangle\right)\right].$$ Since $I^2$ is generated in equal degrees $4$ and the analytic spread $\ell_R(I^2)= \ell_S(I^2S) = 3$, so $e_3\left(\mathfrak{n} \vert \langle \left(I^{2}S\right)_{4}\rangle\right) = 0$. From Macaulay2 computations we get that $e_3\left(\mathfrak{n} \vert \langle \left(\widetilde{I^{2}}S\right)_{4}\rangle\right) = 4$. Thus, $\varepsilon(I)=1/2$ and it agrees with value obtained in \cite[Example $2.4$]{DC4}.
\end{example}

The situation is very simple when $R$ is standard graded polynomial ring in two variables.

\begin{remark}[With hypotheses as in Setup \ref{setup}]\label{Tess2}
Assume that $R$ is a unique factorization domain (UFD). Let $I$ be a height one homogeneous ideal in $R$ which is not principal. Since $R$ is a UFD, we can write $I = (f)\cdot Q$, where $f$ is the greatest common divisor of all nonzero elements in $I$ and $Q$ is a proper homogeneous ideal of height at least two. Further assume that $Q$ is $\mathfrak{m}$-primary (which always happen when $\dim R=2$). Then for all $v\geq 0$, $$(f)^v \subseteq \widetilde{I^v} = \left((f)^v\cdot Q^v\right) \colon_R \mathfrak{m}^{\infty} \subseteq (f)^v \colon_R \mathfrak{m}^{\infty} = (f)^v,$$ where the last equality holds as $R$ is a UFD. Note that the saturated Rees algebra $\mathcal{S}(I)$ is standard graded. Moreover, $$\varepsilon(I) = \lim_{v\to\infty}\dfrac{\lambda_R\left((f)^v/(f)^vQ^v\right)}{v^2/2} = \lim_{v\to\infty}\dfrac{\lambda_R\left(R/Q^v\right)}{v^2/2} = e(Q).$$
\end{remark}

\begin{example}
Let $k$ be a field and consider the monomial ideal $I=(X^2Y,XY^3,Y^5)$ in the polynomial ring $R=k[X,Y]$ with $\mathfrak{m}=(X,Y)$. Then for all $v>0$, $$I^v = (Y)^v\cdot \left(X^2, XY^2, Y^4\right)^v.$$ Let $Q$ denote the $\mathfrak{m}$-primary component $\left(X^2, XY^2, Y^4\right)$ of $I$ and observe that $(X^2,Y^4) \subseteq Q$ is a reduction of $Q$. As $(X^2,Y^4)$ is a parameter ideal in $R$, we have $e(Q) = e\left((X^2,Y^4)\right) = \lambda_R\left(R/(X^2,Y^4)\right)=8$. Thus, $\varepsilon(I) = e(Q) = 8$ by Remark \ref{Tess2}.
\end{example}

\begin{example}\label{grifo}
Let $k$ be an algebraically closed field such that $\mathrm{char}\;k\neq 2$. Let $R=k[X,Y,Z]/(XY-Z^2)$ with $\mathfrak{m}=(x,y,z)$. Consider the prime ideal $P=(x,z) \subseteq R$. The Jacobian criterion shows that the singular locus of $\mathrm{Spec}\;R$ is $\{\mathfrak{m}\}$. So $R$ satisfies $R_1$. Since $R$ is also Cohen-Macaulay, it follows that $R$ is a normal domain.

We now claim that $\widetilde{P^{2v}}=(x^v)$ for all $v>0$. Clearly, $z^{2v}=x^vy^v \in P^{2v} \subseteq P^{(2v)}=\widetilde{P^{2v}}$. Since $y^v\not\in P$ therefore $x^v\in P^{(2v)}$. For the converse, first observe that $x,y$ is an $R$-regular sequence hence $x^a,y^b$ is also an $R$-regular sequence for all positive integers $a, b$. For every integer $v>0$, $$P^{2v} = (x^2, xz, z^2)^v = (x^2, xz, xy)^v \subseteq (x^v) \implies \widetilde{P^{2v}} \subseteq (x^v)\colon_R \mathfrak{m}^{\infty} = (x^v)\colon_R (y)^{\infty} = (x^v).$$ This proves our claim. We shall now apply Theorem \ref{epNoeth} with $v_0=2$ and $\beta_{v_0} = 2$. So, $$\varepsilon(P) = \dfrac{1}{2^2}\cdot \left[e_2\left(\mathfrak{n} \vert \langle \left(\widetilde{P^{2}}S\right)_{2}\rangle\right) - e_2\left(\mathfrak{n} \vert \langle \left(P^{2}S\right)_{2}\rangle\right)\right].$$ Macaulay2 computations show that $e_2\left(\mathfrak{n} \vert \langle \left(\widetilde{P^{2}}S\right)_{2}\rangle\right) = 2$. As $P^2$ is generated in equal degrees $2$ and the analytic spread $\ell_R(P^2)= \ell_S(P^2S) = 2$, so $e_2\left(\mathfrak{n} \vert \langle \left(P^{2}S\right)_{2}\rangle\right) = 0$. Thus, $\varepsilon(P)=1/2$.
\end{example}

\begin{remark}
In Lemma \ref{noetheps}, suppose that $\alpha$ and $\alpha^{\prime}$ are the maximum generating degrees of $I$ and $I'$ respectively. Clearly, $\max\{\alpha, \alpha'\} \leq \beta$. From \eqref{amaorel}, we get for $d=2$,
\[\lim_{v\to\infty}\dfrac{\lambda_R\left({I^{\prime}}^v/I^v\right)}{v^2/2} = e_2(\mathfrak{n}| \langle J'_{\alpha'}\rangle)-e_2(\mathfrak{n}| \langle J_{\alpha}\rangle)-2(\alpha'-\alpha)e_1(\mathfrak{n}| \langle J_{\alpha}\rangle)-(\alpha'-\alpha)^2e_0(\mathfrak{n}| \langle J_{\alpha}\rangle).\]
This description could be more useful for computational purposes, since the number of generators of $J_{\alpha}$ (resp. $J'_{\alpha'}$) is less than that of $J_{\beta}$ (resp. $J'_{\beta}$). One can now update Theorem \ref{epNoeth} accordingly and apply it on Example \ref{grifo}. Using Macaulay2, we get
\[e_2\left(\mathfrak{n} \vert \widetilde{P^{2}}S\right)=0,~ e_2\left(\mathfrak{n} \vert P^{2}S\right)=0,~ e_1\left(\mathfrak{n} \vert P^{2}S\right)=2, ~\mbox{and}~ e_0\left(\mathfrak{n} \vert P^{2}S\right)=2.\]
%\[\frac{1}{4} \cdot \big[0-0-2(1-2)2 -(1-2)^2 2\big]=\frac{1}{2}.\]
Thus the new description again gives $\varepsilon(I)=1/2$.
\end{remark}

\section{Epsilon multiplicity in graded dimension two}\label{epdim2}

\s{\bf Setup}.\label{hyp}
Suppose that $R=\oplus_{u\geq 0}R_u$ is a standard graded finitely generated algebra over a field $R_0=k$ with homogeneous maximal ideal $\mathfrak{m}=\oplus_{u\geq 1}R_u$. Assume that $R$ is a domain and $\dim R\geq 2$. Let $I\subseteq R$ be a nonzero homogeneous ideal such that $I$ is not $\mathfrak{m}$-primary. Assume that $I$ is minimally generated by homogeneous elements $f_1,\ldots,f_s$ of degrees $b_1,\ldots,b_s$ respectively with $b_1\leq \cdots \leq b_s$.

\begin{lemma}[With hypotheses as in Setup \ref{hyp}]\label{analyticallyirred}
Then the local ring $R_{\mathfrak{m}}$ is analytically irreducible.
\end{lemma}
\begin{proof}
Let $\widehat{R}$ denote $\mathfrak{m}$-adic completion of $R$. From \cite[Proposition $10.22 (ii)$]{atiyah} it follows that that there exists a graded isomorphism $$\mathrm{gr}_{\mathfrak{m}}\left(R\right) \cong \mathrm{gr}_{\mathfrak{m}\widehat{R}}\left(\widehat{R}\right)$$ of the associated graded rings. Since $R$ is standard graded we also have an isomorphism $$\mathrm{gr}_{\mathfrak{m}}\left(R\right) \cong R$$ of graded rings. Now observe that $\mathrm{gr}_{\mathfrak{m}\widehat{R}}(\widehat{R}) \cong R$ is a domain and $\bigcap\limits_{v\geq 0} \mathfrak{m}^v\widehat{R} = (0)$ since $\widehat{R}$ is complete with respect to the $\mathfrak{m}$-adic topology. Therefore, $\widehat{R}$ is a domain by \cite[Chapter VIII, Theorem $1$, Page $249$]{Zar}.
\end{proof}

\begin{lemma}[With hypotheses as in Setup \ref{hyp}]\label{linear}
Then there exists an integer $h>0$ such that $$\widetilde{I^{hv}} \subseteq I^v$$ for all integers $v\geq 0$. In particular, $$(\widetilde{I^v})_u = 0 \;\;\text{for all integers} \;\; u<\left\lfloor \frac{v}{h}\right\rfloor b_1 \;\;\text{and}\;\; v> 0.$$
\end{lemma}
\begin{proof}
We first claim that $$\bigcap_{v\geq 0}\left(I^v\widehat{R} \colon_{\widehat{R}}\; {\mathfrak{m}\widehat{R}}^{\infty}\right) = (0).$$ Let us assume that this statement is false. Then there exists a nonzero element $f \in \bigcap\limits_{v\geq 0}\left(I^v\widehat{R} \colon_{\widehat{R}}\; {\mathfrak{m}\widehat{R}}^{\infty}\right)$. By a result of Katz and McAdam \cite[Remark $1.5$]{KM} there exists an integer $c>0$ such that $$I^v\widehat{R} \colon_{\widehat{R}}\; {\mathfrak{m}\widehat{R}}^{\infty} = I^v\widehat{R} \colon_{\widehat{R}}\; \mathfrak{m}^{cv}\widehat{R}$$ for all integers $v\geq 0$. This implies that $f\cdot \mathfrak{m}^{cv}\widehat{R} \subseteq I^v\widehat{R}$ for all $v\geq 0$. Since $\widehat{R}$ is a domain by Lemma \ref{analyticallyirred} and $f\neq 0$, we conclude from \cite[Corollary $6.8.12$]{HS} that $$\mathfrak{m}^{c}\widehat{R} \subseteq \overline{I\widehat{R}}.$$ In other words, $\mathrm{rad}\left(I\widehat{R}\right) = \mathrm{rad}\left(\overline{I\widehat{R}}\right) = \mathfrak{m}\widehat{R}$ which contradicts our hypothesis that $\mathrm{rad}(I)\neq \mathfrak{m}$. Therefore, our claim must be correct. The first assertion of our lemma is now a direct consequence of Swanson's theorem \cite[Theorem $3.2$]{IS2}. It also implies the second assertion of the lemma.
\end{proof}

\begin{lemma}[With hypotheses as in Setup \ref{hyp}]\label{regul}
Then there exists an integer $u_0\geq 0$ such that $$(\widetilde{I^v})_u = (I^v)_u$$ for all integers $u\geq b_sv + u_0$ and $v\geq 0$.
\end{lemma}
\begin{proof}
For every integer $v\geq 0$, consider the short exact sequence $$0 \longrightarrow I^v \longrightarrow  \widetilde{I^v} \longrightarrow \widetilde{I^v}/I^v \longrightarrow 0$$ of graded $R$-modules. Now from the induced long exact sequence of graded local cohomology modules, we obtain that $$0 \longrightarrow \widetilde{I^v}/I^v \longrightarrow H^1_{\mathfrak{m}}\left(I^v\right) \longrightarrow H^1_{\mathfrak{m}}(\widetilde{I^v}) \longrightarrow 0$$ is exact and for every integer $i\geq 2$ there are graded isomorphisms $$H^i_{\mathfrak{m}}\left(I^v\right) \cong H^i_{\mathfrak{m}}(\widetilde{I^v}).$$ We also know from \cite[Theorem $3.2$]{TW} that exists an integer $u_0\geq 0$ such that $$\max\left\{j \mid \left(H^i_{\mathfrak{m}}\left(I^v\right)\right)_{j-i}\neq 0\;\text{and}\;i\geq 0\right\}=:\mathrm{reg}(I^v) \leq b_sv+u_0$$ for integers $v\geq 0$. Hence, if $u\geq b_sv+u_0$ then $\left(H^1_{\mathfrak{m}}\left(I^v\right)\right)_u = 0$ which implies that $(\widetilde{I^v}/I^v)_u = 0$. This proves our claim.
\end{proof}

\s {\bf Riemann-Roch theorems for integral projective curves}.
Let $X$ be an integral projective curve (possibly singular) over an algebraically closed field $k$. For a coherent sheaf $\mathcal{F}$ of $\mathcal{O}_X$-modules, define $h^i(X,\mathcal{F}) = \dim_k H^i\left(X,\mathcal{F}\right)$. The Euler characteristic $\chi(\mathcal{F})$ of $\mathcal{F}$ is defined as $$\chi(\mathcal{F}) = h^0\left(X,\mathcal{F}\right) - h^1\left(X,\mathcal{F}\right).$$ If $D$ is a Cartier divisor on $X$ then $\chi(\mathcal{O}_X(nD))$ is a linear polynomial in $n$ with integer coefficients, see \cite[Chapter VIII, Proposition $1.2$]{altman}. The \emph{degree} of $D$, denoted by $\deg D$, is defined as the leading coefficient of the polynomial $\chi(\mathcal{O}_X(nD))$. The degree of a Cartier divisor depends only on its linear equivalence class. Moreover, if $D_1$ and $D_2$ are two Cartier divisors on $X$ then $$\deg (D_1\pm D_2) = \deg D_1 \pm \deg D_2,$$ see \cite[Chapter VIII, Proposition $1.5$]{altman}. Consequently, the degree function induces a homomorphism $$\deg \colon \mathrm{CaCl}\;X \longrightarrow \mathbb{Z}$$ of abelian groups, where $\mathrm{CaCl}\;X$ is the group of Cartier divisor classes modulo principal divisors.

\begin{theorem}[Riemann-Roch]\label{RR}
Let $D$ be a Cartier divisor on $X$. Then $$\chi(\mathcal{O}_X(D)) = \deg D + \chi(\mathcal{O}_X) = \deg D + 1 - g,$$ where $g=h^1(X,\mathcal{O}_X)$. Moreover, there exists a canonical coherent $\mathcal{O}_X$-module $\omega_X$ such that the vector spaces $$H^1(X,\mathcal{O}_X(D)) \quad \text{and}\quad H^0(X,\omega_X\otimes \mathcal{O}_X(-D))$$ and the vector spaces $$H^0(X,\mathcal{O}_X(D)) \quad \text{and}\quad H^1(X,\omega_X\otimes \mathcal{O}_X(-D))$$ are canonically dual to each other. In particular, $\chi(\mathcal{O}_X(D)) = \chi(\omega_X\otimes \mathcal{O}_X(-D))$.
\end{theorem}
\begin{proof}
The proofs of the above statements can be found in Serre's book \cite[Chapter IV]{serre}.
\end{proof}

If $D$ is an effective Cartier divisor on $X$ then $\deg D = h^0(X,\mathcal{O}_D)$. To see this, consider the short exact sequence $$0 \longrightarrow \mathcal{O}_X(-D) \longrightarrow \mathcal{O}_X \longrightarrow \mathcal{O}_D \longrightarrow 0$$ and use the additivity of the Euler characteristic $\chi(-)$ on short exact sequences and the Riemann-Roch formula.

\begin{lemma}\label{vanishing}
 Let $D$ be a Cartier divisor on $X$. Then the following statements are true:
 \begin{enumerate}
  \item[$(a)$] If $h^0(X,\mathcal{O}_X(D))> 0$ then $\deg D\geq 0$.
  \item[$(b)$] If $\deg D > 2g-2$ then $h^1(X,\mathcal{O}_X(D))=0$.
  \item[$(c)$] Let $\left\{D_n\right\}_{n\in \mathbb{N}}$ be a collection of Cartier divisors on $X$ such that $\left\{\deg D_n\right\}_{n\in \mathbb{N}}$ form a bounded sequence. Then the sequence $\left\{h^0(X,D_n)\right\}_{n\in \mathbb{N}}$ is also bounded.
 \end{enumerate}
\end{lemma}
\begin{proof}
\begin{itemize}
 \item[$(a)$] If $h^0(X,\mathcal{O}_X(D))>0$ then $D$ is linearly equivalent to some effective divisor $E$. Thus $$\deg D = \deg E = h^0(X,\mathcal{O}_E)\geq 0.$$
 \item[$(b)$] This is proven in \cite[Lecture $11$, Page $80$]{mumford} by using the duality statement of the Riemann-Roch theorem.
 \item[$(c)$] Since the sequence $\left\{\deg D_n\right\}_{n\in \mathbb{N}}$ is bounded, we may choose an effective divisor $E$ such that $\deg D_n + \deg E>2g-2$ for all $n\geq 0$. Let $\alpha = \limsup_{n\in \mathbb{N}} \deg D_n$. Now use part $(b)$ and the Riemann-Roch formula to conclude that $$h^0(X,\mathcal{O}_X(D_n)) \leq h^0(X,\mathcal{O}_X(D_n+E)) = \deg D_n + \deg E + 1-g \leq \alpha + \deg E + 1-g.$$
\end{itemize}
\end{proof}

The following lemma is crucial for our main result.

\begin{lemma}[With hypotheses as in Setup \ref{hyp}]\label{castel}
Further assume that $k$ is an algebraically closed field and $\dim R = 2$. Then there exist positive integers $\eta_0, \eta_1$ and nonnegative integers $g, v_0$ such that for all integers $u\geq 0$ and $v\geq v_0$,
\begin{equation*}
 \dim_k (\widetilde{I^v})_u = \begin{cases}
0 & \mathrm{if}\;\;u<\frac{\eta_0}{\eta_1}v,\\
\eta_1u - \eta_0v +1-g  & \mathrm{if}\;\;u>\frac{\eta_0}{\eta_1}v + \frac{2g-2}{\eta_1}
\end{cases}.
\end{equation*}
Moreover, $\eta_0 = -e_0(R[It])$ and $\eta_1 = e_1(R[It])$.
\end{lemma}
\begin{proof}
Let $V = \mathrm{Proj}~R$ be the corresponding projective curve with a very ample invertible sheaf $\mathcal{O}_V(1)$. Let $\mathcal{I}$ be the ideal sheaf associated to the ideal $I$ on $V$. For every integer $v\geq 0$, there is a graded exact sequence of $R$-modules $$0 \longrightarrow I^v \longrightarrow \bigoplus_{u=0}^{\infty}H^0\left(V,\mathcal{I}^v\otimes \mathcal{O}_V(u)\right) \longrightarrow H^1_{\mathfrak{m}}\left(I^v\right) \longrightarrow 0.$$ Moreover, the graded short exact sequence $$0 \longrightarrow I^v \longrightarrow R \longrightarrow R/I^v \longrightarrow 0$$ induces a graded long exact sequence $$0 \longrightarrow H^0_{\mathfrak{m}}\left(R/I^v\right) \longrightarrow H^1_{\mathfrak{m}}\left(I^v\right) \longrightarrow H^1_{\mathfrak{m}}\left(R\right) \longrightarrow \cdots$$ of graded local cohomology modules. There exists an integer $u_0\geq 0$ (for instance, $u_0=\mathrm{reg}(R)$) such that $\left(H^1_{\mathfrak{m}}\left(R\right)\right)_u = 0$ for all $u\geq u_0$. As a consequence, we obtain that
\begin{equation}\label{mainineq0}
 H^0\left(V,\mathcal{I}^v\otimes \mathcal{O}_V(u)\right) = (\widetilde{I^v})_u
\end{equation}
for all $u\geq u_0$ and $v\geq 0$. Let $$\pi\colon X = \mathbf{Proj}\left(\oplus_{v\geq 0}\mathcal{I}^vt^v\right) \longrightarrow V$$ be the blowing up of $V$ along $\mathcal{I}$ with exceptional divisor $E$. Let $H$ be the pullback of a hyperplane section on $V$. Therefore, $$\mathcal{O}_X(H) = \pi^*\mathcal{O}_V(1)\quad \text{and} \quad\mathcal{O}_X(-E) = \mathcal{I}\mathcal{O}_X.$$ From \cite[Lemma $5.4.24$]{Laz} we know that there exists an integer $v_0\geq 0$ such that $$\pi_*\mathcal{O}_X(-vE) = \mathcal{I}^{v}$$ for all $v\geq v_0$. From the projection formula and equation \eqref{mainineq0}, we then see that
\begin{equation}\label{mainineq(-1)}
 H^0\left(X,\mathcal{O}_X(uH-vE)\right) = H^0\left(V,\mathcal{I}^v\otimes \mathcal{O}_V(u)\right) = (\widetilde{I^v})_u
\end{equation}
for all $u\geq u_0$ and $v\geq v_0$. Note that $X$ is an integral projective curve over $k$. The Riemann-Roch formula (see Theorem \ref{RR}) combined with Lemma \ref{vanishing} $(a)$ and $(b)$ give us that
\begin{equation}\label{mainineq(-2)}
 h^0\left(X,\mathcal{O}_X\left(uH-vE\right)\right) = \begin{cases}
0 & \mathrm{if}\;\;u<\frac{v\deg E}{\deg H}\\
u\deg H-v\deg E+1-g  & \mathrm{if}\;\;u>\frac{v\deg E}{\deg H} + \frac{2g-2}{\deg H}
\end{cases},
\end{equation}
where $g=h^1(X,\mathcal{O}_X)$. In view of \eqref{mainineq(-2)} and Lemma \ref{linear}, we may increase $v_0$ and achieve \eqref{mainineq(-1)} for all $u\geq 0$ and $v\geq v_0$. We know from Lemma \ref{regul} and Theorem \ref{hoang-trung} that there exists an integer $u_0^{\prime}\geq 0$ such that for all integers $u\geq b_sv+u_0^{\prime}$ and $v\geq 0$, $$\dim_k (\widetilde{I^v})_u = \dim_k (I^v)_u = e_1(R[It])u + e_0(R[It])v + \eta,$$ where $\eta \in \mathbb{Z}$ is a constant. By comparing the above expression with \eqref{mainineq(-2)}, we obtain that $$\deg H = e_1(R[It]), \quad \deg E = -e_0(R[It]) \quad \text{and} \quad 1-g = \eta.$$ The assertions of our lemma now follow \eqref{mainineq(-2)} and \eqref{mainineq(-1)}.
\end{proof}

\begin{remark}
Let $I$ and $R$ be as in Setup \ref{hyp}. If $\ell_R(I)<\dim R$ then the $\varepsilon$-function $\lambda_R\left(H^0_{\mathfrak{m}}\left(R/I^v\right)\right)$ is eventually given by a polynomial of degree at most $\ell_R(I)-1$, see \cite[Theorem $2.1$]{DC4}.
\end{remark}

We now aim to prove our main result. Its proof is inspired by \cite[Theorem $2.3$]{DC3}.

\begin{theorem}[With hypotheses as in Setup \ref{hyp}]\label{rationalep}
Further assume that $k$ is an algebraically closed field and $\dim R = 2$. Then the following statements hold:
\begin{enumerate}
 \item[$(a)$] If the ideal $I$ has analytic spread two, then for all integers $v\gg 0$, $$\lambda_R\left(H^0_{\mathfrak{m}}\left(R/I^v\right)\right) = a_0v^2 + a_1(v)v + a_2(v),$$ where $a_0\in \mathbb{Q}_{>0}$ is a constant, $a_1 \colon \mathbb{N} \to \mathbb{Q}$ is a periodic function and $a_2 \colon \mathbb{N} \to \mathbb{Q}$ is a bounded function.
 \item[$(b)$] Let $S = R[U]$, where $U$ is a variable with $\deg U=1$, $\mathfrak{n} = \mathfrak{m} + (U)$ is the homogeneous maximal ideal of $S$ and $J=IS$ is the extended ideal in $S$. Then $$\varepsilon(I) = \dfrac{\left(e_1\left(S[Jt]\right)\right)^2}{e_2\left(S[Jt]\right)} - e_0\left(S[Jt]\right) = \dfrac{\left(e_1\left(\mathfrak{n}\vert\left<J_{b_s}\right>\right)\right)^2}{e_0\left(\mathfrak{n}\vert\left<J_{b_s}\right>\right)} - e_2\left(\mathfrak{n}\vert\left<J_{b_s}\right>\right).$$
\end{enumerate}
\end{theorem}

\begin{proof}
By Lemma \ref{regul} there exists an integer $c>b_s$ such that $$(\widetilde{I^v})_u = \left(I^v\right)_u$$ for all integers $v\geq 0$ and $u\geq cv$. Define $$\sigma(v) = \sum_{u=0}^{cv}\dim_{k}(\widetilde{I^v})_u  \quad\text{and}\quad \tau(v) = \sum_{u=0}^{cv}\dim_{k}\left(I^v\right)_u.$$ This allows us to write
\begin{equation*}
 \lambda_R\left(H^0_{\mathfrak{m}}\left(R/I^v\right)\right) = \sigma(v) - \tau(v).
\end{equation*}
Recall the proof of Lemma \ref{castel} and we know that there exists an integer $v_0\geq 0$ such that $$H^0\left(X,\mathcal{O}_X(uH-vE)\right) = (\widetilde{I^v})_u$$ for all $u\geq 0$ and $v\geq v_0$. Consider the function $$\delta(v) = \sum_{u=\left\lceil \frac{v\deg E}{\deg H}\right\rceil}^{\left\lceil \frac{v\deg E}{\deg H} + \frac{2g-1}{\deg H}\right\rceil -1}h^0\left(X,\mathcal{O}_X\left(uH-vE\right)\right).$$ Lemma \ref{vanishing} $(c)$ shows that $\delta \colon \mathbb{N} \to \mathbb{N}$ is a bounded function. Now for all $v\geq v_0$, one has
\begin{align}\label{mainineq2}
 \sigma(v) = \sum_{u=0}^{cv}h^0\left(X,\mathcal{O}_X\left(uH-vE\right)\right)\nonumber &= \delta(u) + \sum_{u=\left\lceil \frac{v\deg E}{\deg H} + \frac{2g-1}{\deg H}\right\rceil}^{cv}\left(u\deg H-v\deg E+1-g\right)\nonumber\\
 &= \left[\dfrac{c^2\deg H}{2} - c\deg E + \dfrac{(\deg E)^2}{2\deg H}\right]v^2 + \gamma_1(v)v + \gamma_2(v),
\end{align}
where $\gamma_1 \colon \mathbb{N} \to \mathbb{Q}$ is a periodic function and $\gamma_2 \colon \mathbb{N} \to \mathbb{Q}$ is a bounded function. Since $c>b_s$, it follows from an argument used in the proof of \cite[Theorem $2.3$]{DC3} that $$\tau(v) = \sum_{u=0}^{cv}\dim_k \left(I^v\right)_u = \dim_k \left(J^v\right)_{cv}$$ for all $v\geq 0$. From \eqref{HT1} it immediately follows that $\tau(v)$ eventually agrees with a quadratic polynomial in $v$. In particular, by possibly increasing $v_0$, we have that for all $v\geq v_0$,
\begin{equation}\label{mainineq3}
 \tau(v) = \left[\dfrac{c^2e_2\left(S[Jt]\right)}{2} + ce_1\left(S[Jt]\right) + \dfrac{e_0\left(S[Jt]\right)}{2}\right]v^2 + \text{lower degree terms in $v$}.
\end{equation}
For proving part $(b)$ of our theorem, we observe from \eqref{mainineq2} and \eqref{mainineq3} that $$\varepsilon(I) = \lim_{v\to\infty}\dfrac{\sigma(v)-\tau(v)}{v^2/2} = \left(\deg H - e_2\left(S[Jt]\right)\right)c^2 - 2\left(\deg E + e_1\left(S[Jt]\right)\right)c + \dfrac{(\deg E)^2}{\deg H} - e_0\left(S[Jt]\right).$$ The above expression is true for all sufficiently large values of $c$ and so we must have $$\deg H = e_2\left(S[Jt]\right) \; \text{and} \; \deg E = -e_1\left(S[Jt]\right) \;\text{which implies}\; \varepsilon(I) = \dfrac{\left(e_1\left(S[Jt]\right)\right)^2}{e_2\left(S[Jt]\right)} - e_0\left(S[Jt]\right).$$ From Proposition \ref{mixedrel0} and Lemma \ref{mixedrel} we obtain the following set of identities:
\begin{align}\label{relate}
 e_2(S[Jt]) &= e_0\left(\mathfrak{n} \vert \langle J_{b_s}\rangle\right),\\
 e_1(S[Jt]) &= e_1\left(\mathfrak{n} \vert \langle J_{b_s}\rangle\right) - b_se_0\left(\mathfrak{n} \vert \langle J_{b_s}\rangle\right),\nonumber\\
 e_0(S[Jt]) &= e_2\left(\mathfrak{n} \vert \langle J_{b_s}\rangle\right) - 2b_se_1\left(\mathfrak{n} \vert \langle J_{b_s}\rangle\right) + b_s^2 e_0\left(\mathfrak{n} \vert \langle J_{b_s}\rangle\right).\nonumber
\end{align}
Using these relations, we further conclude that $$\varepsilon(I) = \dfrac{\left(e_1\left(\mathfrak{n}\vert\left<J_{b_s}\right>\right)\right)^2}{e_0\left(\mathfrak{n}\vert\left<J_{b_s}\right>\right)} - e_2\left(\mathfrak{n}\vert\left<J_{b_s}\right>\right).$$ This proves part $(b)$ of our theorem. Now, part $(a)$ follows by subtracting \eqref{mainineq3} from \eqref{mainineq2}.
\end{proof}
%
% \begin{remark}
% Recall the proof of Lemma \ref{castel}. Let $\varphi \colon X^{\prime} \longrightarrow V$ denote the normalization of $V$. Any ideal sheaf on $X^{\prime}$ is invertible as $X^{\prime}$ is a nonsingular projective curve. In particular, $\mathcal{I}\mathcal{O}_{X^{\prime}}$ is an invertible sheaf on $X^{\prime}$. Note that $\varphi \colon X^{\prime} \to V$ can also be realized as blowing up of some ideal sheaf on $V$. It then follows from the universal property of blowing up \cite[Chapter II, Proposition $7.14$]{hartshorne} that there exists a unique morphism $\psi \colon X^{\prime} \longrightarrow X$ such that $\varphi = \psi \circ \pi$. In fact, $\psi \colon X^{\prime} \to X$ is also the normalization of $X$. Let $E^{\prime} = \psi^*E$ and $H^{\prime} = \psi^*H$. Then $$\psi^*\mathcal{O}_X(H) = \mathcal{O}_{X^{\prime}}(H^{\prime}) = \varphi^*\mathcal{O}_V(1) \quad \text{and} \quad \psi^*\mathcal{O}_X(-E) = \mathcal{O}_{X^{\prime}}(-E^{\prime}) = \mathcal{I}\mathcal{O}_{X^{\prime}}$$ and since $\psi$ is birational, so $$\deg H^{\prime} = \deg H = e_2\left(S[Jt]\right) \quad \text{and} \quad \deg E^{\prime} = \deg E = -e_1\left(S[Jt]\right).$$ Here, $e_2\left(S[Jt]\right)$ is equal to the Hilbert-Samuel multiplicity of $R_{\mathfrak{m}}$. Therefore, $$\varepsilon(I) = \dfrac{\left(e_1(S[Jt]\right)^2}{e_2(S[Jt])} - e_0\left(S[Jt]\right) = \dfrac{(\deg E^{\prime})^2}{\deg H^{\prime}} - e_0\left(S[Jt]\right).$$
% \end{remark}

\begin{corollary}[With hypotheses as in Theorem \ref{rationalep}]\label{equaldeg}
Further assume that the ideal $I$ is generated in a single degree $b>0$. Then $$\varepsilon(I) = \dfrac{\left(e_1\left(\mathfrak{m}\vert I\right)\right)^2}{e_0\left(\mathfrak{m}\vert I\right)}.$$
\end{corollary}
\begin{proof}
Here, we use the formula obtained in part $(b)$ of Theorem \ref{rationalep}. Since the ideal $I$ is generated in equal degrees so $\left<J_b\right> = J$. As $\ell_R(I) = \ell_S(J)\leq 2$, so $e_2\left(\mathfrak{n}\vert J\right) = 0$ by \cite[Theorem $2.7$]{katz-verma}. Moreover, from the arguments used in the proof of Corollary \ref{superficial0}, we see that $$e_0\left(\mathfrak{n}\vert J\right) = e_0\left(\mathfrak{m}\vert I\right)\;\text{and}\;e_1\left(\mathfrak{n}\vert J\right) = e_1\left(\mathfrak{m}\vert I\right).$$ This implies our claim.
\end{proof}

\begin{remark}[With hypotheses as in Theorem \ref{rationalep}]\label{epscomp}
Suppose that $P$ is a homogeneous prime ideal of height one in $R$. By the projective Nullstellensatz, $P$ corresponds to a closed point $p$ on the projective curve $V=\mathrm{Proj}\;R$ and $P$ is generated by linear forms in $R$. Let $\mathcal{O}_{V,p}$ denote the local ring of $V$ at $p$, and there is an isomorphism $$\mathcal{O}_{V,p} \cong R_{(P)},$$ where $R_{(P)}$ is the homogeneous localization of $R$ at $P$. The integral closure $\overline{\mathcal{O}}_{V,p}$ of $\mathcal{O}_{V,p}$ is a Dedekind domain. So there exists a unique factorization $$P\overline{\mathcal{O}}_{V,p} = Q_1^{\alpha_1}\cdots Q_r^{\alpha_r},$$ where $Q_1,\ldots,Q_r$ are the distinct maximal ideals in $\overline{\mathcal{O}}_{V,p}$ lying over $P\mathcal{O}_{V,p}$ and $\alpha_1,\ldots,\alpha_r \in \mathbb{Z}_{\geq 1}$, see \cite[Corollary $9.4$]{atiyah}. Now, it follows that $$e\left(\mathcal{O}_{V,p}\right) = \sum_{i=1}^r \alpha_i,$$ see \cite[Theorem $11.2.7$]{HS}. From the proofs of Lemma \ref{castel} and Theorem \ref{rationalep} we conclude that $$e_0\left(\mathfrak{m}\vert P\right) = e(R) \;\;\text{and} \;\; e_1\left(\mathfrak{m}\vert P\right) = e(R) - e\left(\mathcal{O}_{V,p}\right).$$ In particular, Corollary \ref{equaldeg} implies that $$\varepsilon_R(P) = \dfrac{\left(e(R)-e\left(\mathcal{O}_{V,p}\right)\right)^2}{e(R)}.$$ This expression will help one understand the collection of pairs $(R,P)$ where $R$ satisfies the hypotheses of Theorem \ref{rationalep} and $P\subseteq R$ a height one homogeneous prime ideal, for which $\varepsilon_R(P)$ is constant. For instance, consider the collection of all such pairs $(R,P)$ where $R$ is normal and $e(R)$ is constant, say $e$, then $$\varepsilon_R(P) = \dfrac{\left(e-1\right)^2}{e} = e-2+\dfrac{1}{e}.$$
\end{remark}

\begin{corollary}[With hypotheses as in Theorem \ref{rationalep}]\label{periodic1}
Further assume that $k$ is an algebraically closed field of characteristic zero. If the ideal $I$ has analytic spread two, then for all $v\gg 0$, $$\lambda_R\left(H^0_{\mathfrak{m}}\left(R/I^v\right)\right) = a_0v^2 + a_1(v)v + a_2(v),$$ where $a_0\in \mathbb{Q}_{>0}$ is a constant and $a_i \colon \mathbb{N} \to \mathbb{Q}$ is a periodic function for $i=1,2$.
\end{corollary}
\begin{proof}
Recall the proofs of Lemma \ref{castel} and Theorem \ref{rationalep}. It was shown in \eqref{mainineq(-1)} that there exists an integer $v_0\geq 0$ such that for all integers $u\geq 0$ and $v\geq v_0$, $$H^0\left(X,\mathcal{O}_X(uH-vE)\right) = H^0\left(V,\mathcal{I}^v\otimes \mathcal{O}_V(u)\right) = (\widetilde{I^v})_u.$$ A result of Cutkosky (refer to the comment after \cite[Theorem $4.1$]{DC7}) shows that there exists a rational expression of the bivariate Poincar\'{e} series
\begin{equation}\label{mainineq4}
\sum\limits_{(u,v)\in\NN^2}h^0\left(X,\mathcal{O}_X\left(uH-vE\right)\right)t_1^{u}t_2^{v} = \dfrac{P_1(t_1,t_2)}{\prod_{i=1}^{r_1} \left(1-t_1^{a_{i1}}t_2^{a_{i2}}\right)},
\end{equation}
where $P_1(t_1,t_2) \in \mathbb{Q}[t_1,t_2]$ and $(a_{i1},a_{i2})\in\NN^2\setminus\{(0,0)\}$ for all $i=1,\ldots,r_1$. In view of \eqref{mainineq(-1)}, the Hilbert series of the bigraded saturated Rees algebra $\oplus_{(u,v)\in\NN^2}(\widetilde{I^v})_u t_1^{u}t_2^{v}$, i.e.,
\begin{align*}
 \sum_{v=0}^{\infty}\sum_{u=0}^{\infty} \dim_k (\widetilde{I^v})_u t_1^{u}t_2^{v} =  \sum_{v=0}^{\infty}\sum_{u=0}^{\infty}h^0\left(X,\mathcal{O}_X\left(uH-vE\right)\right)t_1^{u}t_2^{v} &- \sum_{v=0}^{v_0-1}\sum_{u=0}^{\infty}h^0\left(X,\mathcal{O}_X\left(uH-vE\right)\right)t_1^{u}t_2^{v}\\ &\qquad + \sum_{v=0}^{v_0-1}\sum_{u=0}^{\infty}\dim_k (\widetilde{I^v})_u t_1^{u}t_2^{v}
\end{align*}
also has a rational expression similar to the one given in \eqref{mainineq4}. Also, the Hilbert series of the bigraded Rees algebra $\oplus_{(u,v)\in\NN^2}\left(I^v\right)_u t_1^{u}t_2^{v}$ has a rational expression of the form
\begin{equation*}
 \sum\limits_{(u,v)\in\NN^2}\dim_k \left(I^v\right)_u t_1^{u}t_2^{v} = \dfrac{P_2(t_1,t_2)}{\prod_{i=1}^{r_2} \left(1-t_1^{c_i}t_2\right)},
\end{equation*}
where $P_2(t_1,t_2) \in \mathbb{Q}[t_1,t_2]$ and $c_1,\ldots,c_{r_2}$ are some positive integers. Then notice that
\begin{align*}
 \sum\limits_{(u,v)\in\NN^2} \dim_k \left(H^0_{\mathfrak{m}}\left(R/I^v\right)\right)_u t_1^{u}t_2^{v} &= \sum\limits_{(u,v)\in\NN^2} \dim_k (\widetilde{I^v})_u t_1^{u}t_2^{v} - \sum\limits_{(u,v)\in\NN^2}\dim_k \left(I^v\right)_u t_1^{u}t_2^{v}
\end{align*}
again has a rational expression similar to the one given in \eqref{mainineq4}. For all $v\geq 0$, the expression $\sum_{u=0}^{\infty} \dim_k \left(H^0_{\mathfrak{m}}\left(R/I^v\right)\right)_u t_1^{u}t_2^{v}$ is a finite sum. So, the reduced rational expression of the series $\sum_{(u,v)\in\NN^2} \dim_k \left(H^0_{\mathfrak{m}}\left(R/I^v\right)\right)_u t_1^{u}t_2^{v}$ cannot have a pole at $t_1=1$. Putting $t_1=1$ gives an expression $$\sum\limits_{v=0}^{\infty} \lambda_R \left(H^0_{\mathfrak{m}}\left(R/I^v\right)\right) t_2^{v} = \sum\limits_{v=0}^{\infty} \left(\sum_{u=0}^{\infty}\dim_k \left(H^0_{\mathfrak{m}}\left(R/I^v\right)\right)_u\right)t_2^v = \dfrac{P_3(t_2)}{\prod_{i=1}^{r_3} \left(1-t_2^{e_i}\right)},$$ where $P_3(t_2) \in \mathbb{Q}[t_2]$ and $e_1,\ldots,e_{r_3}$ are some positive integers. In other words, $\lambda_R \left(H^0_{\mathfrak{m}}\left(R/I^v\right)\right)$ eventually agrees with a numerical quasi-polynomial in $v$ of degree two.
\end{proof}

\begin{remark}
The saturated Rees algebra $\mathcal{S}(I)=\oplus_{v\geq 0}\widetilde{I^v}t^v$ of $I$ may be non-Noetherian even in the situation considered in Corollary \ref{periodic1}. To see this, take $R$ to be the homogeneous coordinate ring of an elliptic curve, say $R=\frac{\mathbb{C}[X,Y,Z]}{(X^3+Y^3+Z^3)}$, and let $P$ be a height one homogeneous prime ideal in $R$ which has infinite order in the divisor class group, see \cite{rees}. Then the algebra $\mathcal{S}(P)$ is non-Noetherian.
\end{remark}

\begin{remark}
The bounded function $a_2 \colon \mathbb{N} \to \mathbb{Q}$,  introduced in Theorem \ref{rationalep}$(a)$, need not be periodic when $\mathrm{char}\;k>0$ and $k$ has positive transcendence degree over its prime subfield. One may use the counterexample to Zariski's Riemann-Roch problem in positive characteristic, due to S. D. Cutkosky and V. Srinivas \cite[Section $6$]{CS}, to construct such an example where the function $a_2(v)$ is bounded but is not eventually periodic.
\end{remark}

\s {\bf Some examples}.
Throughout this section, $\mathbb{C}$ will denote the field of complex numbers.
\begin{example}\cite[Exercise $3.7$]{fulton}
Let $n\geq 1$ be an odd integer and consider the polar equation $r=\sin(n \theta)$. We use the trigonometric identity,

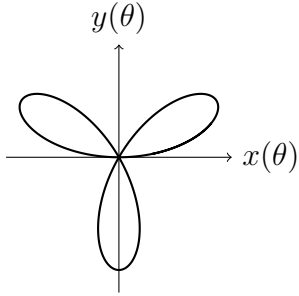
\begin{minipage}{.27\columnwidth}
\begin{tikzpicture}[scale=1.5]
\draw[->] (-1,0) -- (1,0) node[right] {$x(\theta)$};
\draw[->] (0,-1.2) -- (0,1) node[above] {$y(\theta)$};
\draw[thick,variable=\t,domain=0.1:3.5,samples=500]
plot ({sin((3*\t) r)*cos((\t) r)},{sin((3*\t) r)*sin((\t) r)});
\draw (0,-1.5) node {\footnotesize$\bigl(x(\theta),y(\theta)\bigr) = $};
\draw (0,-1.8) node {\footnotesize$(\sin(3 \theta)\cos \theta, \sin(3\theta)\sin \theta)$};
\end{tikzpicture}
\end{minipage}
\begin{minipage}{.7 \columnwidth}
\[\sin (n \theta)=\sum_{t=0}^n \binom{n}{t} (\cos \theta)^t (\sin \theta)^{n-t} \sin \frac{(n-t)\pi}{2},\] to convert this polar equation to its Cartesian form by setting $X=r \cos \theta$ and $Y=r \sin \theta$.
We get, \[(X^2+Y^2)^{\frac{n+1}{2}}=\sum_{t=0}^n \binom{n}{t} X^t Y^{n-t} \sin \frac{(n-t)\pi}{2}.\]
Now homogenize this affine equation and consider the irreducible homogeneous polynomial
\[f_n(X,Y,Z) := (X^2+Y^2)^{\frac{n+1}{2}} -\sum_{t=0}^n \binom{n}{t} \sin \frac{t \pi}{2} X^{n-t}Y^tZ\]
\end{minipage}
in $\mathbb{C}[X,Y,Z]$. Let $R=\mathbb{C}[X,Y,Z]/(f_n(X,Y,Z))$ and observe that $V = \mathrm{Proj}\; R$ is a projective plane curve of degree $n+1$ with an ordinary $n$-fold point at $p=[0\colon 0\colon 1]$. Let $P=(x,y)$ be the homogeneous prime ideal in $R$ corresponding to the point $p$. We shall now use the formula given in Remark \ref{epscomp} to compute the $\varepsilon$-multiplicity of $P$. Clearly, $e(R) = n+1$ and $e\left(R_{(P)}\right) = n$ and so $$\varepsilon(P) = \dfrac{\left(e(R)-e\left(R_{(P)}\right)\right)^2}{e(R)} = \dfrac{(n+1-n)^2}{n+1} = \dfrac{1}{n+1}.$$
\end{example}

\begin{example}\label{exdense}
Fix two integers $p$ and $q$ such that $0<p<q$ and $\gcd(p,q)=1$. Let $R = \mathbb{C}[X,Y,Z]/(X^q-Y^{q-p}Z^p)$. Note that $V= \mathrm{Proj}\;R$ is a projective plane curve of degree $q$ and
$$\mathrm{Sing}(V) = \begin{cases}
                      \emptyset &\text{if}\;p=1\;\text{and}\;q=2\\
                      \left\{[0\colon 0 \colon 1]\right\} &\text{if}\;p=1\;\text{and}\;q\geq 3\\
                      \left\{[0\colon 1 \colon 0]\right\} &\text{if}\;p\geq 2\;\text{and}\;q=p+1\\
                      \left\{[0\colon 1 \colon 0],[0\colon 0 \colon 1]\right\} &\text{if}\;p\geq 2\;\text{and}\;q\geq p+2
                     \end{cases}.
$$
Let $P=(x,y)$ be the homogeneous prime ideal in $R$ corresponding to the point $[0\colon 0\colon 1]$ on $V$. Consider the open affine subset $$D_{+}(z) = \mathrm{Spec}\;\dfrac{\mathbb{C}[X,Y]}{(X^q-Y^{q-p})} \cong \mathrm{Spec}\;\mathbb{C}[s^{q-p},s^q]$$ containing the point $[0\colon 0\colon 1]$. We have $e(R) = q$ and $e\left(R_{(P)}\right) = e\left(\mathbb{C}[s^{q-p},s^q]_{(s^{q-p},s^q)}\right) = q-p$. Now use the formula in Remark \ref{epscomp} to get $$\varepsilon(P) = \dfrac{\left(e(R)-e\left(R_{(P)}\right)\right)^2}{e(R)} = \dfrac{p^2}{q}.$$
\end{example}

\begin{example}\cite{Tei86}
Let $$R=\dfrac{\mathbb{C}[X,Y,Z,W]}{\left(Y^2W-X^3, Z^2W^4-X^5Y\right)}.$$ Let $V = \mathrm{Proj}\;R \subseteq \mathbb{P}^3_{\mathbb{C}}$ be the associated projective curve and observe that $V$ is singular at the point $p=[0\colon 0\colon 0\colon 1]$. Let $P=(x,y,z)$ be the homogeneous prime ideal in $R$ corresponding to the point $p \in V$. Consider the open affine subset $$D_{+}(w) = \mathrm{Spec}\;\dfrac{\mathbb{C}[X,Y,Z]}{\left(Y^2-X^3, Z^2-X^5Y\right)} \cong \mathrm{Spec}\;\mathbb{C}[s^4,s^6,s^{13}]$$ containing $p$. The integral closure of $\mathbb{C}[s^4,s^6,s^{13}]$ in its field of fractions $\mathbb{C}(s)$ is $\mathbb{C}[s]$. Thus, $e\left(R_{(P)}\right) = e\left(\mathbb{C}[s^4, s^6, s^{13}]_{(s^4,s^6,s^{13})}\right) =  4$, by \cite[Theorem $11.2.7$]{HS}.

For convenience, let us denote by $R^{\prime}=\mathbb{C}[X,Y,Z,W]$ and set $f=Y^2W-X^3$ and $g=Z^2W^4-X^5Y$. Note that $f,g$ is an $R^{\prime}$-regular sequence with $\deg f=3$ and $\deg g=6$. One can verify that the Hilbert series of $R$ is
\[H_R(t) = \dfrac{(1-t^6)(1-t^3)}{(1-t)^4} = \dfrac{(1+t+\cdots+t^5)(1+t+t^2)}{(1-t)^2}.\]
Thus, $e(R)=(1+t+\cdots+t^5)(1+t+t^2)\big|_{t=1}=18$. Now by Remark \ref{epscomp}, $$\varepsilon(P) = \dfrac{\left(e(R)-e\left(R_{(P)}\right)\right)^2}{e(R)} = \dfrac{98}{9}.$$
\end{example}

\begin{example}\cite[Example $2.3$]{GS21}
Consider the elliptic curve $V$ with homogeneous coordinate ring $R=\frac{\mathbb{C}[X,Y,Z]}{\left(X^3-Y^2Z-2Z^3\right)}$. Consider the prime ideal $P=(x-3z,y-5z)\subseteq R$ which corresponds to the closed point $p=[3\colon 5\colon 1]\in V$. In this case, the saturated Rees algebra $\mathcal{S}(P)$ of $P$ is non-Noetherian since $p$ has infinite order with respect to the group law on $V$, see \cite{rees}.

By contrast, consider the point $q=[2\colon 3\colon 1]$ on the elliptic curve with homogeneous coordinate ring $R^{\prime}=\frac{\mathbb{C}[X,Y,Z]}{\left(X^3-Y^2Z+Z^3\right)}$. The point $q$ has has order six with respect to the group law of this curve. The defining ideal of the point $q$ is the prime ideal $Q=(x-2z,y-3z)\subseteq R^{\prime}$. Macaulay2 computations give $$\widetilde{Q^6}=(12x^2-6xy+y^2-6xz-6yz+9z^2)$$ which is a principal ideal. Moreover, $\widetilde{Q^{6v}} = \left(\widetilde{Q^6}\right)^v$ for all $v\geq 0$. This shows that the saturated Rees algebra $\mathcal{S}(Q)$ is Noetherian.

We can apply Remark \ref{epscomp} to either of the two situations considered above and get $$\varepsilon_R(P) = \varepsilon_{R^{\prime}}(Q) = \dfrac{(3-1)^2}{3} = \dfrac{4}{3}.$$
\end{example}

\s {\bf Discussions}.

Let $R$ be as in Theorem \ref{rationalep}. In Remark \ref{epscomp}, we gave a simple description for $\varepsilon(I)$ when $I$ is a height one homogeneous prime ideal in $R$. It is now natural to ask whether one can extract similar formula for $\varepsilon(I)$ when $I=\cap_{i=1}^rP_i^{a_i}$, where $P_i$'s are some height one homogeneous prime ideals in $R$ and $a_i$'s are some positive integers. From Theorem \ref{rationalep}, we have $$e_2(S[Jt]) = \deg H = e(R)\quad \text{and}\quad e_1(S[Jt]) = -\deg E = -\sum\limits_{i=1}^r a_i\cdot e\left(R_{(P_i)}\right).$$ Notice that if $I$ is equigenerated then $e_0(S[Jt])$ can be explained by means of \eqref{relate}, where $e_2(\mathfrak{n}\vert\langle J_{b_s}\rangle)=e_2(\mathfrak{n}|J)=0$. Otherwise, if $I$ is generated in different degrees, $e_2(S[Jt])$ will depend on the configuration of the closed points $p_i$ (associated to $P_i$'s) in $\mathrm{Proj}\;R$. Therefore, a theoretical understanding of $e_2(S[Jt])$ seems difficult.

We now describe the epsilon multiplicity of certain fat points. Some statements are purely based on computations in Macaulay2 and we don't have any concrete proof right now.

Consider the ring $R=\mathbb{C}[X,Y,Z]/(X^3+Y^3+Z^3)$ and pick the ideals $P_1=(x+y,z)$ and $P_2=(x+z,y)$ associated to the points $p_1=[1:-1:0]$ and $p_2=[1:0:-1]$ respectively on the elliptic curve $\mathrm{Proj}\;R$. We would like to compute the $\varepsilon$-multiplicities of the ideal $I=P_1^a\cap P_2^b$ for different values of $a$ and $b$. By symmetry, we can assume that $1 \leq a \leq b$. Note that \[\deg H = e_0(\mathfrak{n}\vert J_{b_s})=e(R)=3 \quad \mbox{and}\quad \deg E=a+b.\] The following data is obtained using Macaulay2:

\begin{center}
\begin{tabular}{|p{2cm}|p{2cm}|p{2cm}|}
\hline
$(a,b)$ & $b_s$ & $e_2(\mathfrak{n}\vert \langle J_{b_s}\rangle)$\\
\hline
$(1,1)$ & $2$ & $4$ \\
$(1,2)$ & $3$ & $8$ \\
$(1,3)$ & $4$ & $12$ \\
$(2,2)$ & $3$ & $5$ \\
$(2,3)$ & $4$ & $9$ \\
$(3,5)$ & $6$ & $14$ \\
\hline
\end{tabular}
\end{center}

This information raises the following question:

\begin{question}
Is $b_s=b+1$ and $e_2(\mathfrak{n}\vert J_{b_s})=4b-3a+3$?
\end{question}

For now, assume that this question has an affirmative answer. Then from \eqref{relate}, we conclude $$e_1(\mathfrak{n}\vert J_{b_s})=b_s\deg H-\deg E=2b-a+3.$$ Moreover, by Theorem \ref{rationalep} $$\varepsilon(I)=\frac{1}{3}(2b-a)^2+a.$$

\vspace{0.15cm}
Next, we consider the ring $R=\mathbb{C}[X,Y,Z]/(X^3-Y^2Z)$. Choose ideals $P_1=(x,y)$ and $P_2=(x,z)$ corresponding to the points $p_1=[0:0:1]$ and $p_2=[0:1:0]$, respectively. Let $I=P_1^a\cap P_2^b$. It is easy to verify that $$\deg H = e_0(\mathfrak{n}\vert J_{b_s})=e(R)=3\quad \mbox{and}\quad \deg E = a\cdot e(R_{(P_1)}) + b\cdot e(R_{(P_2)}) = 2a+b.$$ The following data is obtained from Macaulay2:

\begin{center}
\begin{tabular}{|p{2cm}|p{2cm}|p{2cm}|}
\hline
$(a,b)$ & $b_s$ & $e_2(\mathfrak{n}\vert \langle J_{b_s}\rangle)$ \\
\hline
$(1,1)$ & $2$ & $3$ \\
$(1,2)$ & $3$ & $7$  \\
$(2,1)$ & $2$ & $0$  \\
$(1,3)$ & $4$ & $11$  \\
$(3,1)$ & $3$ & $0$  \\
$(2,2)$ & $3$ & $3$  \\
$(2,3)$ & $4$ & $7$  \\
$(3,2)$ & $3$ & $0$  \\
$(2,5)$ & $6$ & $15$  \\
$(5,2)$ & $5$ & $0$  \\
$(3,3)$ & $4$ & $3$  \\
$(3,5)$ & $6$ & $11$  \\
$(5,3)$ & $5$ & $0$  \\
$(5,7)$ & $8$ & $11$  \\
\hline
\end{tabular}
\end{center}

The above outputs suggest the following:
\begin{question}
Is $b_s=b+1$ and $e_2(\mathfrak{n}|J_{b_s})=4b-4a+3$ if $1 \leq a\leq  b$? Again, if $1 \leq b <a$ then is it true that $I$ generated in equal degree $b$?
\end{question}
 Notice in the second case $e_2(\mathfrak{n}\vert \langle J_{b_s}\rangle)=0$. Provided we have an affirmative answer to the above question, by \eqref{relate} and Theorem \ref{rationalep} it follows that
\[\begin{cases}
e_1(\mathfrak{m}|I)=3a-(2a+b)=a-b ~\mbox{ and }~ \varepsilon(I)=\frac{(a-b)^2}{3} & \quad \mbox{ if } 1 \leq b <a\\
e_1(\mathfrak{n}|J_{b_s})=3(b+1)-(2a+b)=2b-2a+3 ~\mbox{ and } ~\varepsilon(I)=\frac{4(b-a)^2}{3} & \quad \mbox{ if } 1 \leq a\leq  b.
  \end{cases}
\]

% This computation can also be done using simple algebraic arguments which we discuss below.
% \begin{claim*}
% For all integers $n\geq 1$, we have $\widetilde{I^{qn}} = \left(y^{(q-p)n}\right) = \left(\widetilde{I^q}\right)^n$.
% \end{claim*}
%
% We now parametrize $R$ via the map $$\theta \colon \mathbb{C}[X,Y,Z] \longrightarrow \mathbb{C}[U,V]$$ given by $\theta(X) = U^{q-p}V^p$, $\theta(Y) = U^q$ and $\theta(Z) = V^q$. Observe that $\mathrm{ker}\;\theta = (X^q-Y^{q-p}Z^p)$ and $\theta$ induces a graded isomorphism $$R = \dfrac{\mathbb{C}[X,Y,Z]}{(X^q-Y^{q-p}Z^p)} \cong \mathbb{C}[U^q, U^{q-p}V^p, V^q].$$ Let $\Gamma$ denote the affine semigroup in $\mathbb{Z}^2$ generated by the vectors $\{(q,0), (q-p,p), (0,q)\}$. Then $R$ is isomorphic to the semigroup ring $\mathbb{C}[\Gamma] = \mathbb{C}\left[U^{\alpha}V^{\beta} \mid (\alpha,\beta)\in \Gamma\right]$. Let $\overline{\Gamma}$ denote the saturation of the semigroup $\Gamma$, i.e., $\overline{\Gamma} = \left(\mathbb{R}_{\geq 0}\cdot \Gamma\right) \cap \mathbb{Z}^2$. Then the integral closure $\overline{R}$ of $R$ is isomorphic to $\mathbb{C}[\overline{\Gamma}] = \mathbb{C}\left[U^{\alpha}V^{\beta} \mid (\alpha,\beta)\in \overline{\Gamma}\right]$. Also, for all $n\geq 1$, $$\lambda_R\left(H^0_{\mathfrak{m}}\left(R/I^{qn}\right)\right) = \lambda_R\left(\widetilde{I^{qn}}/I^{qn}\right) = \lambda_{\mathbb{C}[\Gamma]}\left(\dfrac{\left(U^{q(q-p)}\right)^{n}}{\left(U^{q-p}V^p, U^q\right)^{qn}}\right)$$

\section{An algorithm for computing \texorpdfstring{$\varepsilon$}{epsilon}-multiplicities}\label{M2_algorithm}

Let $R$ be a finitely generated graded algebra over a field $k$, that is, $R=k[x_1, \ldots, x_n] \cong k[X_1, \ldots, X_n]/J$ for some homogeneous ideal $J$ in $k[X_1, \ldots, X_n]$. Suppose that $\deg x_i=d_i$ for $i=1, \ldots, n$. Let $I=(g_1, \ldots, g_s)$ be a homogeneous ideal in $R$ with $\deg g_i=e_i$ for $i=1, \ldots, s$. In this section, we restrict ourselves to the case, where the saturated Rees algebra $\mathcal{S}(I)$ is Noetherian.

We first recollect some examples. If $I$ is an edge ideal, then some bounds on the generating degrees of $\mathcal{S}(I)$ are known from \cite{HHT}. If $I$ is the defining ideal of the space monomial curve $k[s^a, s^b, s^c]$ then Herzog and Ulrich \cite{HU90} developed a theory on self-linked curve singularities and gave a criterion for the symbolic Rees algebra $\mathcal{R}_s(I)$ to be generated in degree two. This was succeeded by Goto, Nishida and Shimoda \cite{GNS91}, where they found conditions under which $\mathcal{R}_s(I)$ is generated in degrees at most three. Some extensions of this theory to the case when the generating degrees of $\mathcal{R}_s(I)$ is at most four, can be found in \cite{Reed09}. The \emph{generation type} of a symbolic Rees algebra $\mathcal{R}_s(I)$ is defined as \[\mathrm{gt}\big(\mathcal{R}_s(I)\big):=\inf\left\{r\in\mathbb{N} \mid \mathcal{R}_s(I)=R[It, I^{(2)}t^2, \ldots, I^{(r)}t^r]\right\}\] and the \emph{standard Veronese degree} of an ideal $I$ is \[\mathrm{svd}(I):=\inf\left\{r\in\mathbb{N} \mid \big(I^{(r)}\big)^n=I^{(rn)} \;\mbox{for all}\; n \geq 1\right\}.\]
Readers may refer to \cite{GS21} for a short survey regarding certain bounds on these invariants.

\vspace{0.15cm}
Suppose that $\mathcal{S}(I)$ is generated in degrees at most $s$. For each $i=1,\ldots,s$, write $\widetilde{I^i}=\langle f_{ij} \mid j=1,\ldots,m_i \rangle$ for some homogeneous elements $f_{ij}\in R$. Consider the set of indeterminates $\underline{Y} = \{Y_{ij} \mid i=1,\ldots,s, j=1,\ldots,m_i \}$ and $\underline{X}=\{X_1, \ldots, X_n\}$. Define an $R$-algebra homomorphism
\[\varphi\colon k[\underline{X}][\underline{Y}] \longrightarrow \mathcal{S}(I) \subseteq R[T]\] given by $\varphi(Y_{ij}) = f_{ij}T^i$ for $i=1,\ldots,s$, $j=1,\ldots,m_i$ and $\varphi(X_t) = x_t$ for $t=1, \ldots, n.$ Then $\mathcal{S}(I) \cong k[\underline{X},\underline{Y}] / \mathrm{ker}(\varphi).$ There is a natural bigraded structure on $R[T]$ where $\deg x_t=(0,d_t)$ and $\deg T=(1,0)$. We further assign a bigraded structure on $S:=k[\underline{X},\underline{Y}]$ by setting $\deg X_t=(0,d_t)$ and $\deg Y_{ij}=(i,\deg f_{ij})$. Then, $\varphi$ is a bigraded ring homomorphism. Thus, $\mathrm{ker}\; \varphi$ has an induced bigraded structure. Similarly, one can get a bigraded ring isomorphim $\mathcal{R}(I) \cong k[\underline{X}][Z_1, \ldots, Z_s]/\mathcal{L}$ by sending $Z_i \to g_iT$, where $\mathcal{R}(I)$ denotes the Rees algebra of $I$.

This allows us to compute their bivariate Hilbert series $$H_{\mathcal{S}(I)}(t_0, t_1) := \sum_{(u,v)\in\mathbb{N}^2}\dim_k \left(\widetilde{I^v}\right)_ut_0^vt_1^u\quad \text{and}\quad H_{\mathcal{R}(I)}(t_0, t_1) := \sum_{(u,v)\in\mathbb{N}^2}\dim_k \left(I^v\right)_ut_0^vt_1^u,$$ using Macaulay2. Observe that \[\left(\mathcal{S}(I)/\mathcal{R}(I)\right)_{(u,v)}=\left(\widetilde{I^v}/I^v\right)_u = \left(H^0_{\mathfrak{m}} \left(R/I^v\right)\right)_u\] for every $(u,v) \in \mathbb{N}^2$. Let $W(u,v)=\dim_k \mathcal{S}(I)_{(u,v)}-\dim_k \mathcal{R}(I)_{(u,v)}$. Then \[\sum_{(u,v) \in \mathbb{N}^2} W(u,v) t_0^vt_1^u = H_{\mathcal{S}(I)}(t_0, t_1) - H_{\mathcal{R}(I)}(t_0, t_1) =: \mathcal{H}(t_0, t_1).\] As $\lambda_R\left(H^0_{\mathfrak{m}}(R/I^v)\right)< \infty$, so \[\mathcal{H}(t_0,1)=\sum_{v \geq 0}\left(\sum_{u\geq 0} W(u,v)\right) t_0^v = \sum_{v\geq 0} \lambda_R\left(H^0_{\mathfrak{m}} \left(R/I^v\right)\right)t_0^v.\]

We now write a script in Macaulay2 that we use to compute the $\varepsilon$-multiplicity and $\varepsilon$-function.

\noindent
{\bf Input:} A sequence $W$ of an ideal $I$ in a finitely generated graded (not necessarily standard graded) algebra $R$ over a field $k$ and a positive integer $n$ such that $\mathcal{S}(I)$ is generated in degrees at most $n$.

\noindent
{\bf Output:} The rational function associated to the series $\sum_{v \geq 0} \lambda_R(H^0_{\mathfrak{m}}(R/I^v)) t^v$, where $t$ is an auxiliary variable and $\mathfrak{m}$ is the unique homogeneous maximal ideal of $R$.

\vspace{0.15cm}
\noindent
{\bf Script:}
{\footnotesize
\begin{multicols}{2}
\begin{verbatim}
--W=(I,n)
epsilonSeries=W->(
n := W#1;
I := W#0;
-- to check if $I$ is an ideal
if class I =!= Ideal then error "The
sequence is incorrect";
R := ring W#0;
-- to check if the analytic spread of $I$
is equal to the dimension of $R$
if  analyticSpread I =!= dim R then error "The
analytic spread of the ideal is not maximum";
m := ideal vars R;
k := coefficientRing R;
degR := flatten degrees ideal vars R;
bidegR := flatten toList
                  (0..#degR-1)/(i->{0,degR#i});
D := join(bidegR, {{1,0}});
J = (i)->(saturate(I^i, m));
S := k[gens R, t, Degrees=>D];
MJ = (i) ->(mingens (sub(J(i),S)*t^i));
DJ = (i) ->(degrees ideal MJ(i));
MR := matrix{delete(t, gens S)};
a = MR; for j from 1 to n do
			  (b = a|MJ(j); a = b;);
M := a;
NV :=sum toList apply(1..n, i->#DJ(i));
ND := join apply(1..n, i->DJ(i));
DD:= join(bidegR, ND);
U := k[gens R, Y_1..Y_NV, Degrees=> DD];
KS := ker map(S, U, M);
HKS := reduceHilbert hilbertSeries KS;
KRs := reesIdeal I;
HKRs := reduceHilbert hilbertSeries KRs;
A := QQ[T_0, T_1];
HNS := sub((toList HKS)#0, A);
HDS := sub(value (toList HKS)#1, A);
HKS1 := HNS/HDS;
HNRs := sub((toList HKRs)#0, A);
HDRs := sub(value (toList HKRs)#1, A);
HKRs1 := HNRs/HDRs;
H := HKS1 - HKRs1;
HD := sub(denominator H, A/ideal(T_1-1));
HN := sub(numerator H, A/ideal(T_1-1));
B := ZZ[T_0];
HD1 := sub(HD, B);
HN1 := sub(HN, B);
G := gcd (HD1, HN1);
HD2 := sub(HD1/G, B);
HDF := factor HD2;
HN2 := sub(HN1/G, B);
HNF := factor HN2;
H1 := HNF/HDF;
H1)
\end{verbatim}
\end{multicols}
}

\begin{remark}\label{epscompute}
We may use a number of softwares to compute the partial fraction decomposition of $H1:=\mathcal{H}(t_0, 1)$. One can verify that if `$c$' is the largest integer appearing in the expression as a power of $\frac{1}{(1-t)}$ then $c=d+1$ where $d:=\dim R$. Moreover, $\varepsilon(I)$ is the coefficient of $\frac{1}{(1-t)^c}$. Whenever, the above script produces an output, one can compute the epsilon multiplicity $\varepsilon(I)$ along with the epsilon function $\varepsilon_I(n)$.
\end{remark}

\begin{example}
Let $k$ be a field and $R=k[X,Y,Z]$ is a graded polynomial ring with $\deg X=3, \deg Y=4,$ and $\deg Z=5$. Consider, the ring homomorphism $\varphi\colon k[X,Y,Z] \to k[s]$ given by $\varphi(X)=s^3$, $\varphi(Y)=s^4$ and $\varphi(Z)=s^5$. Let $P=\mathrm{ker}\;\varphi$ and it can be shown that $P=I_2\left(\begin{bmatrix}
X^2 & Y & Z\\
Z & X & Y
\end{bmatrix}\right)$, see \cite[Page $137$]{Kun85}. Then $k[s^3, s^4, s^5] \cong R/P$. Let $m=(x,y,z)$ be the unique homogeneous maximal ideal of $R$. By \cite[Corollary 2.12]{HU90}, the symbolic Rees algebra $\mathcal{R}_s(P)$ of the monomial space curve is of the form $\mathcal{R}_s(P)=R\left[Pt, P^{(2)}t^2\right].$
We run the following session in Macaulay2.
{\footnotesize
\begin{multicols}{2}
\begin{verbatim}
i1 : R=QQ[x,y,z, Degrees=>{{3},{4},{5}}];
i2 : M=matrix{{x^2,y,z},{z,x,y}};
i3 : I=minors(2,M);
i4 : isHomogeneous I
o4 = true
i5 : epsilonSeries (I,2)
\end{verbatim}
\end{multicols}
}

\noindent
Given a space curve, one can extract the corresponding matrix using the following commands.
{\footnotesize
	\begin{multicols}{2}
		\begin{verbatim}
	i1 : S=QQ[t];
	i2 : R=QQ[x,y,z, Degrees=>{{3},{4},{5}}];
	i3 : P= ker map(S, R, matrix{{t^3, t^4, t^5}});
	i4 : C=res P;
	i5 : M=transpose C.dd_2;
	     --the corresponding matrix
		\end{verbatim}
	\end{multicols}
}
Here $M=\begin{bmatrix}
-Z & -Y & X\\
-X^2 & -Z & Y
\end{bmatrix}$.
Using an online partial fraction decomposition calculator in WolframAlpha \cite{WolAlp09}, we get
\begin{align}\label{epsilonlast}
\mathcal{H}(t_0,1)&=\frac{t_0^2}{(1-t_0)^4(1+t_0)}\nonumber\\
&=t_0^2\left(\frac{1}{16(1+t_0)}+\frac{1}{16(1-t_0)}+\frac{1}{8(1-t_0)^2}-\frac{3}{4(1-t_0)^{3}}+\frac{1}{2(1-t_0)^{4}}\right)\nonumber\\
&=\sum_{n=0}^\infty \left[\frac{(-1)^n}{16}+\frac{1}{16}+\frac{1}{8}(n+1)-\frac{3}{4}\binom{n+2}{2}+\frac{1}{2}\binom{n+3}{3}\right]t_0^{n+2}.
\end{align}
Observe that the epsilon function
\[\varepsilon_P(n)=\frac{1}{16}+\frac{(n-1)}{8}+\frac{n(n-1)}{8}+\frac{(n+1)n(n-1)}{12}+\frac{(-1)^{n}}{16}\]
for $n \geq 2$. Hence
\[\varepsilon(P)=\lim\limits_{n \to \infty}\dfrac{\lambda_R\left(H^0_\mathfrak{m}(R/P^n)\right)}{n^3/ 3!}=1/2.\] Here, $\varepsilon(P)$ can also be directly computed from the expression \eqref{epsilonlast} by Remark \ref{epscompute}.
\end{example}

\begin{example}\cite[Example 6.17 (1)]{GNS91}
    Consider the ideal $P=I_2\left(\begin{bmatrix}
	X^2 & Y^2 & Z^3\\
	Y & Z & X^3
	\end{bmatrix}\right)$ in $R=k[X,Y,Z]$ which is the defining ideal of the monomial space curve $k[s^{10},s^{11},s^{13}]$. By \cite[Theorem $6.1$]{GNS91}, $\mathcal{S}(P)=R\left[Pt, P^{(2)}t^2, P^{(3)}t^3\right]$. Using Macaulay2 and an online calculator WolframAlpha, we get
	\begin{align*}
	\mathcal{H}(t_0,1)&=\frac{2t_0^2(1+t_0)}{(1-t_0)^4(t_0^2+t_0+1)}\\
	&=2t_0^2\Bigg[-\frac{1}{9(t_0^2+t_0+1)}+\frac{1}{9(1-t_0)^2}+\frac{1}{3(1-t_0)^3}+\frac{2}{3(1-t_0)^4}\Bigg]\\
	&=-\frac{2t_0^2}{9(t_0^2+t_0+1)}+\sum_{n=0}^\infty\Bigg[ \frac{1}{9}(n+1)+\frac{1}{3}\binom{n+2}{2}+\frac{2}{3}\binom{n+3}{3}\Bigg]t_0^{n+2}.
	\end{align*}
	Hence $\varepsilon(I)=2/3$ by Remark \ref{epscompute}.
	\end{example}

\begin{example}
Consider the ideal $P=I_2\left(\begin{bmatrix}
X^2 & Y^2 & Z^3\\
Y & Z^2 & X^2
\end{bmatrix}\right)$ in $R=k[X,Y,Z]$ which is the defining ideal of the monomial space curve $k[s^{11},s^{14},s^{10}]$. In light of \cite[Section $2$]{GNS91} (see also \cite[Notation $2.4$]{Reed09}), we have $\{(\alpha, \beta, \gamma), (\alpha', \beta',\gamma')=\{(2,1,2),(2,2,3)\}\}$. Thus the matrix is of type $I$, see \cite[Page $101$]{GNS91}. Using Macaulay2 we get
	\[\left(Pt\right)^4+Pt \cdot P^{(3)}t^3+ \left(P^{(2)}t^2\right)^2 \subsetneq P^{(4)}t^4.\]
	Thus $\mathcal{S}(P)\neq R\left[Pt, P^{(2)}t^2, P^{(3)}t^3\right]$, see \cite[Proposition $3.1$(3)]{GNS91}. Hence by \cite[Theorem $4.1(h)$]{Reed09}, $\mathcal{S}(P)=R\left[Pt, P^{(2)}t^2, P^{(3)}t^3, P^{(4)}t^4\right]$.
	Using Macaulay2 and an online calculator WolframAlpha, we get
	\begin{align*}
	\mathcal{H}(t_0,1)&=\frac{2t_0^2(2+t_0+2t_0^2)}{(1-t_0)^4(1+t_0^2)(1+t_0)}\\
	&=2t_0^2\Bigg[-\frac{(1+t_0)}{8(1+t_0^2)}-\frac{1}{32(1-t_0)}+\frac{3}{32(1+t_0)}+\frac{3}{16(1-t_0)^2}+\frac{5}{8(1-t_0)^3}+\frac{5}{4(1-t_0)^4}\Bigg].
	\end{align*}
	Hence $\varepsilon(I)=5/4$ by Remark \ref{epscompute}.
\end{example}

\begin{example}\label{I_4_3}\cite[Example $3.9$]{GS21}
	Let $R=k[X_1,\cdots, X_4]$ and consider the monomial ideal
	\begin{align*}
	I:=I_{4,3}&=(X_1,X_2,X_3) \cap(X_1,X_2,X_4)\cap(X_1,X_3,X_4) \cap (X_2,X_3,X_4)\\
	&=(X_3X_4,X_2X_4,X_1X_4,X_2X_3,X_1X_3,X_1X_2)
	\end{align*}
	Clearly, $\mathrm{ht}\; I=3$ and by \cite[Lemma $3.2$]{MMV12}, $\ell_R(I) = \mathrm{rank}\;M_I=4$, where
	\[M_I=\begin{bmatrix}
	0&0&1&0&1&1\\
	0&1&0&1&0&1\\
	1&0&0&1&1&0\\
	1&1&1&0&0&0
	\end{bmatrix}.\]
Therefore, the saturation powers and the symbolic powers of $I$ coincide. By \cite[Example $3.9$]{GS21}, $\mathrm{gt}\big(\mathcal{R}_s(I)\big)=3$, i.e., $\mathcal{S}(I)=R[It, \widetilde{I^2}t^2, \widetilde{I^3}t^3]$. Using Macaulay2 and an online calculator WolframAlpha, we get
	\begin{align*}
	\mathcal{H}(t_0,1)&=\frac{t_0^2(t_0^4+3t_0^3+7t_0^2+5t_0+4)}{(1-t_0)^5(1+t_0)(t_0^2+t_0+1)}\\
	&=t_0^2\Bigg[\frac{t_0+2}{27(t_0^2+t_0+1)}+\frac{35}{216(1-t_0)}+\frac{1}{8(1+t_0)}+\frac{1}{4(1-t_0)^2}\\
	&\hspace{5cm}+\frac{7}{18(1-t_0)^3}-\frac{1}{3(1-t_0)^4}+\frac{10}{3(1-t_0)^5}\Bigg].
	\end{align*}
	Hence $\varepsilon(I)=10/3$ by Remark \ref{epscompute}.
\end{example}

\begin{remark}
The relation described in Theorem \ref{epNoeth} is not practically useful for Example \ref{I_4_3}. From \cite[Example $3.9$]{GS21} we get $\mathrm{svd}(I)$ is divisible by $6$. Whereas, by \cite[Proposition $3.1$ and Remark $3.2$]{GS21}, $\mathrm{svd}(I) \leq \mathrm{gt}\big(\mathcal{R}_s(I)\big)\cdot \mathrm{gt}\big(\mathcal{R}_s(I)\big)!=18$. Using Macaulay2, we get $$\mu(I)=6, \mu\big(\widetilde{I^2}\big)=10, \mu\big(\widetilde{I^3}\big)=19, \mu\big(\widetilde{I^6}\big)=460.$$ Again from the computations in Macaulay2, it seems to suggest that $\mathrm{svd}(I)=6$. If this is true then set $\widetilde{I^6}R[X_5]=\mathbb{J}$, which is generated in degrees $\{8,9,10,11,12\}$. As $\mu\big(\langle \mathbb{J}_{12}\rangle\big)=460$ so Macaulay2 mostly crashes while computing the mixed multiplicities even in this case.
\end{remark}

\begin{example}
 Let $k$ be a field and $M=(X_{ij})_{m \times n}$ be a matrix of indeterminates over $k$ with $m \leq n$. Consider the polynomial ring $R=k[X_{ij} \mid 1 \leq i \leq m, 1 \leq j \leq n]$ and the ideal $I_s$ generated by the $s\times s$-minors of $X$. Then it is well-known that
 \begin{equation}\label{symbolicconca}
  \mathcal{S}(I_s)=R[I_st,I_{s+1}t^2,\ldots, I_mt^{m-s+1}],
 \end{equation}
 see \cite[Proposition $4.1$]{Con98}. Now, take $s=2$ and $m=n=3$. Then $\ell_R(I_2)=9=\dim R$ by \cite[Remark $6.5$]{JMV15} and $\mathrm{ht}\;I_2 = 4$. From \cite[Theorem $4.4$]{JMV15} we get that $I^{(n)}=\widetilde{I^n}$ for all $n \geq 1$. Computing in Macaulay2 and using an online calculator WolframAlpha, we get
		\begin{align*}
		\mathcal{H}(t_0,1)&=\frac{t_0^2}{(1-t_0)^{10}(1+t_0)}\\
		&=\dfrac{1}{1024(1+t_0)}+\frac{1}{1024(1-t_0)}+\frac{1}{512(1-t_0)^2}+\frac{1}{256(1-t_0)^3}+\frac{1}{128(1-t_0)^4}+\frac{1}{64(1-t_0)^5}\\
		&\hspace{4cm} +\frac{1}{32(1-t_0)^6}+\frac{1}{16(1-t_0)^7}+\frac{1}{8(1-t_0)^8}-\frac{3}{4(1-t_0)^9}+\frac{1}{2(1-t_0)^{10}}.
		\end{align*}
		Therefore, $\varepsilon(I)=1/2$ by Remark \ref{epscompute}.
\end{example}

\begin{remark}
 The algorithm will be more effective if one uses the description \eqref{symbolicconca} of $\mathcal{S}(I_s)$ to compute its defining ideal and thereby its bivariate Hilbert Series.
\end{remark}

\begin{example}[Fermat ideals]
Consider the ideal $I=\left(X(Y^3-Z^3), Y(Z^3-X^3), Z(X^3-Y^3)\right)$ in the standard graded polynomial ring $R=\mathbb{C}[X,Y,Z]$, where $\mathbb{C}$ denotes the field of complex numbers. It follows from \cite[Proposition $4.1$]{NS16} that $I^{(3v)} = \left(I^{(3)}\right)^v$ for all integers $v>0$. Besides, $\mathrm{ht}(I)=2$ and $\ell_R(I)=3=\dim R$. Thus $I^{(v)} = \widetilde{I^v}$ for all $v>0$. Set $\mathcal{H}_I^{\langle3\rangle}(t_0, t_1):=H_{\mathcal{S}(I)^{\langle 3 \rangle}} (t_0,t_1)-H_{\mathcal{R}(I)^{\langle 3 \rangle}} (t_0,t_1)$. Then
	\[\mathcal{H}_I^{\langle3\rangle}(t_0, 1)= \sum_{v \geq 0} \lambda_R\left(\frac{I^{(3v)}}{I^{3v}}\right)t_0^v = \sum_{v \geq 0} \lambda_R\left(\frac{\left(I^{(3)}\right)^v}{(I^{3})^v}\right)t_0^v.\]
	Using Macaulay2 and an online calculator WolframAlpha, we get
	\[\mathcal{H}_I^{\langle3\rangle}(t_0, 1)=\frac{4t_0^3+82t_0^2+22t_0}{(1-t_0)^4}=-\frac{4}{1-t_0}+\frac{94}{(1-t_0)^2}-\frac{198}{(1-t_0)^3}+\frac{108}{(1-t_0)^4}.\]
	Hence $\varepsilon(I)=\frac{108}{27}=4$ by Remark \ref{epscompute}.
\end{example}

\bibliographystyle{alpha}
\bibliography{Ref}
\end{document}